\documentclass[12pt]{article}

\pdfoutput=1

\usepackage{amsmath}
\usepackage{amssymb}
\usepackage{enumitem}
\usepackage{graphicx}

\def\B{{\boldsymbol B}}
\def\Be{{\boldsymbol B}^{\rm even}}
\def\Bo{{\boldsymbol B}^{\rm odd}}

\def\de{{\delta_m^{\rm even}}}
\def\dod{{\delta_m^{\rm odd}}}
\def\nice{\displaystyle}

\def\End{{\rm End}}

\def\oM{\overline{\mathcal{M}}}
\def\cM{{\mathcal{M}}}

\def\Z{\mathbb{Z}}
\def\C{\mathbb{C}}
\def\Q{\mathbb{Q}}
\def\d{\partial}
\def\qed{{\hfill $\Diamond$}}
\def\b1{{\mathbf 1}}
\def\cS{{\mathcal S}}
\def\hd{{\widehat \partial}}
\def\hdx{{\widehat dx}}
\def\hdy{{\widehat dy}}
\def\cR{{\mathcal R}}
\def\Aut{{\rm Aut}}
\def\F{{\mathsf{F}}}

\def\E{\mathrm{E}}
\def\n{\mathrm{n}}
\def\L{\mathrm{L}}
\def\V{\mathrm{V}}
\def\H{\mathrm{H}}
\def\g{\mathrm{g}}
\def\rarr{\rightarrow}
\def\W{\mathsf{W}}
\def\D{\mathsf{D}}
\def\P{\mathsf{P}}
\def\tP{\widetilde{\mathsf{P}}}
\def\tW{\widetilde{\W}}

\def\xxx{\mathsf{t}}

\usepackage{theorem}

\newtheorem{theorem}{Theorem}

\newtheorem{proposition}{Proposition}[section]
\newtheorem{corollary}[theorem]{Corollary}
\newtheorem{lemma}[proposition]{Lemma}

{\theorembodyfont{\rmfamily}
\newtheorem{definition}[proposition]{Definition}
\newtheorem{example}[proposition]{Example}
\newtheorem{remark}[proposition]{Remark}

}

\title{Relations on $\oM_{g,n}$ via 3-spin structures}
\author{Rahul Pandharipande, Aaron Pixton, Dimitri Zvonkine}
\date{March 2015}

\begin{document}

\maketitle

\vspace{-20pt}

\begin{abstract} 
Witten's class on the moduli space of 3-spin curves
defines a (non-semisimple) cohomological field theory. 
After a canonical modification, we construct an
associated semisimple CohFT with a non-trivial
vanishing property obtained from the homogeneity
of Witten's class. Using the classification of semisimple
CohFTs by Givental-Teleman, we derive two main results.
The first is an explicit formula in the tautological
ring of $\oM_{g,n}$
for Witten's class. The second, using
the vanishing property, is the construction of
relations in the 
tautological ring of $\oM_{g,n}$.

Pixton has previously conjectured a system of tautological
relations on  $\oM_{g,n}$
(which extends the established
Faber-Zagier relations on $\cM_{g}$).  
Our 3-spin construction exactly yields Pixton's conjectured relations.
As the classification of CohFTs is a topological result
depending upon the Madsen-Weiss theorem (Mumford's conjecture),
our construction proves relations in cohomology.
The study of Witten's class and the associated
tautological relations for $r$-spin curves via a parallel
strategy will be taken up in a following paper.
\end{abstract}

\vspace{-20pt}

\setcounter{tocdepth}{1} 
\tableofcontents


\setcounter{section}{-1}
\section{Introduction}

\subsection{Overview}
The study of relations in the cohomology of the moduli space of curves
was initiated by Mumford \cite{Mum} in the 1980s. While several
classical approaches were applied
with success before, 
the subject has developed rapidly in the last two decades via  natural
connections to topological string theory.

A systematic study by Faber and Zagier of the algebra of $\kappa$ classes
on the moduli space $\cM_g$ of
nonsingular genus $g$ curves led to a conjecture in 2000 of 
a concise
set $\mathsf{FZ}$ of $\kappa$ relations.
A proof of the Faber-Zagier conjecture (in Chow) via the geometry of stable quotients
was given in 2010~\cite{PPFZ}.
In 2012, the second author \cite{Pixton}  conjectured 
a set $\mathsf{P}$ of tautological relations for the moduli spaces $\oM_{g,n}$
of stable curves.
The set $\mathsf{P}$ 
recovers  $\mathsf{FZ}$ when restricted to $\cM_g \subset \oM_g$.

Our main result proves the conjectured relations $\mathsf{P}$  
in the cohomology ring $H^*(\oM_{g,n},\mathbb{Q})$.
By restriction, we obtain a second proof of the Faber-Zagier conjecture
in cohomology.
Are there other relations?
The sets $\mathsf{FZ}$ and $\mathsf{P}$ explain {\em all} presently
known tautological relations on $\cM_g$ and $\oM_{g,n}$
respectively. At least in Chow, the sets
$\mathsf{FZ}$ and $\mathsf{P}$ are conjectured 
to be complete in both cases \cite{PPFZ,Pixton}. 

We study here the geometry of 3-spin curves.
Witten's class on the moduli space of 3-spin curves
defines a non-semisimple cohomological field theory. 
After a canonical modification (obtained by moving to a semisimple point of the
associated Frobenius manifold), we construct a
semisimple CohFT with a non-trivial
vanishing property obtained from the homogeneity
of Witten's class. Using the classification of semisimple
CohFTs by Givental-Teleman \cite{Givental, Teleman}, we derive 
an explicit formula in the tautological
ring of $\oM_{g,n}$
for Witten's $3$-spin class {\em and} 
use the vanishing property 
to establish the relation set $\mathsf{P}$.

\subsection{Stable graphs} \label{stgr}
The boundary
strata of the moduli space of curves correspond
to {\em stable graphs} 
$$\Gamma=(\V, \H,\L, \ \mathrm{g}:\V \rarr \Z_{\geq 0},
\ v:\H\rarr \V, 
\ \iota : \H\rarr \H)$$
satisfying the following properties:
\begin{enumerate}
\item[(i)] $\V$ is a vertex set with a genus function $\g:V\to \Z_{\geq 0}$,
\item[(ii)] $\H$ is a half-edge set equipped with a 
vertex assignment $v:H \to V$ and an involution $\iota$,
\item[(iii)] $\E$, the edge set, is defined by the
2-cycles of $\iota$ in $\H$ (self-edges at vertices
are permitted),
\item[(iv)] $\L$, the set of legs, is defined by the fixed points of $\iota$ and endowed with a bijective correspondence with a set of markings,
\item[(v)] the pair $(\V,\E)$ defines a {\em connected} graph,
\item[(vi)] for each vertex $v$, the stability condition holds:
$$2\g(v)-2+ \n(v) >0,$$
where $\n(v)$ is the valence of $\Gamma$ at $v$ including 
both edges and legs.
\end{enumerate}
An automorphism of $\Gamma$ consists of automorphisms
of the sets $\V$ and $\H$ which leave invariant the
structures $\mathrm{g}$, $\iota$, and $v$ (and hence respect $\E$ and $\L$).
Let $\text{Aut}(\Gamma)$ denote the automorphism group of $\Gamma$.

The genus of a stable graph $\Gamma$ is defined by:
$$\g(\Gamma)= \sum_{v\in V} \g(v) + h^1(\Gamma).$$
A boundary stratum of the moduli space $\oM_{g,n}$ 
of Deligne-Mumford stable curves naturally determines
a stable graph of genus $g$ with $n$ legs by considering the dual graph of a generic pointed curve parametrized by the stratum.

To each stable graph $\Gamma$, we associate the moduli space
\begin{equation*}
\oM_\Gamma =\prod_{v\in \V} \oM_{\g(v),\n(v)}.
\end{equation*}
 Let $\pi_v$ denote the projection from $\oM_\Gamma$ to 
$\oM_{\g(v),\n(v)}$ associated to the vertex~$v$.  There is a
canonical
morphism 
\begin{equation}\label{dwwd}
\xi_{\Gamma}: \oM_{\Gamma} \rarr \oM_{g,n}
\end{equation}
 with image{\footnote{
The degree of $\xi_\Gamma$ is $|\text{Aut}(\Gamma)|$.}}
equal to the boundary stratum
associated to the graph $\Gamma$.  To construct $\xi_\Gamma$, 
a family of stable pointed curves over $\oM_\Gamma$ is required.  Such a family
is easily defined 
by attaching the pull-backs of the universal families over each of the 
$\oM_{\g(v),\n(v)}$  along the sections corresponding to half-edges.

\subsection{Strata algebra} \label{Ssec:introdualgraphs}

Let $\Gamma$ be a stable graph.
A {\em basic class} on $\oM_\Gamma$
is defined to be a product of monomials in $\kappa$ classes{\footnote{Our
convention is 
$\kappa_i= \pi_*(\psi_{n+1}^{i+1})\ \in H^{2i}(\oM_{g,n},\mathbb{Q})$ where
$$\pi: \oM_{g,n+1} \rightarrow \oM_{g,n}$$ is the
map forgetting the marking $n+1$. For a review of $\kappa$ and
and cotangent $\psi$ classes, see \cite{GraPan}.}}
 at each 
vertex of the graph and powers of $\psi$ classes at each half-edge (including the legs),
$$
\gamma = \prod_{v\in \mathrm{V}} 
\prod_{i>0}\kappa_i[v]^{x_i[v]} \ \cdot  
\ \prod_{h \in \mathrm{H}} \psi_{h}^{y[h]}\ 
\in H^*(\oM_\Gamma,\mathbb{Q})\ ,
$$
where $\kappa_i[v]$ is the $i^{\rm th}$ kappa class on $\oM_{\mathrm{g}(v),\mathrm{n}(v)}$.
We impose the condition 
$$
\sum_{i>0} i x_i[v] 
+ \sum_{h\in \H[v]} y[h] \leq 
\text{dim}_{\mathbb{C}}\ \oM_{\mathrm{g}(v), \mathrm{n}(v)} =
3\mathrm{g}(v)-3+ \mathrm{n}(v)
$$
at each vertex to avoid  the trivial vanishing of $\gamma$.
Here, 
$\H[v]\subset \H$
is the set of half-edges (including the legs) incident to~$v$.

Consider the $\mathbb{Q}$-vector space $\cS_{g,n}$ whose basis is given by the isomorphism classes of pairs $[\Gamma,\gamma]$, where $\Gamma$ is a stable graph of genus $g$ with $n$
legs and $\gamma$ is a basic class on $\oM_\Gamma$. 
Since there are only finitely many pairs $[\Gamma, \gamma]$ up to isomorphism, $\cS_{g,n}$ is finite dimensional.

A product on $\cS_{g,n}$ is defined by intersection theory
with respect to the morphisms \eqref{dwwd} to $\oM_{g,n}$.
Let 
$$
[\Gamma_1,\gamma_1], \ [\Gamma_2,\gamma_2]\in \cS_{g,n}
$$
be two basis elements. The fiber product of $\xi_{\Gamma_1}$ and $\xi_{\Gamma_2}$ over $\oM_{g,n}$ is canonically described as a disjoint union 
of $\xi_{\Gamma}$ for stable graphs $\Gamma$ endowed with contractions{\footnote{
If there are several different pairs of contractions from a given $\Gamma$, the corresponding $\xi_\Gamma$ appears with multiplicity.}}
 onto $\Gamma_1$ and $\Gamma_2$. More precisely, the set of edges $E$ of $\Gamma$ should be represented as a union of two (not necessarily disjoint) subsets, 
$$E = E_1 \cup E_2,$$ in such a way that
 $\Gamma_1$ is obtained by contracting all the edges outside $E_1$ and 
$\Gamma_2$ is obtained by contracting all edges outside $E_2$ (see Proposition 9 in the Appendix of \cite{GraPan}). 
The intersection of $\xi_{\Gamma_1}$ and $\xi_{\Gamma_2}$ in $\oM_{g,n}$ is then canonically given by Fulton's excess theory as a sum of elements in $\cS_{g,n}$.
 We define
$$
[\Gamma_1,\gamma_1] \cdot [\Gamma_2,\gamma_2] = 
\sum_{\Gamma} [\Gamma, \gamma_1 \gamma_2  \varepsilon_\Gamma]
$$
where 
$$
\varepsilon_\Gamma = \prod_{e \in E_1 \cap E_2} -(\psi_e' + \psi_e'')
$$ 
is the excess class. Here, $\psi'_e$ and $\psi''_e$ are the two 
cotangent line classes corresponding to the two half-edges of the edge $e$.

A case of particular importance for us is when $\Gamma_2$ has only a
single edge. The set $E_2$ must consist of a single element $e$, while $E_1$ may be either $E$ or $E \setminus \{ e \}$. The above product then yields the restriction
of a basic class to a  boundary divisor.

Via the above intersection product, $\cS_{g,n}$ is a finite dimensional $\mathbb{Q}$-algebra, called the {\em strata algebra} \cite{Pixton}.
Push-forward along $\xi_\Gamma$ defines a canonical ring homomorphism 
$$q:\cS_{g,n} \to H^*(\oM_{g,n},\mathbb{Q}), \ \ \
q( [\Gamma,\gamma]) = \xi_{\Gamma*}(\gamma)$$ from the strata algebra to the 
cohomology ring. By definition, the image of $q$ is the tautological ring 
$RH^*(\oM_{g,n})$. An element of the kernel of $q$ is called a {\em tautological relation}.

Each basis element $[\Gamma,\gamma]$ has a degree grading given
by the number of edges of $\Gamma$ plus the usual (complex)
degree of $\gamma$,
$$\text{deg} [\Gamma,\gamma] = | \E| + \text{deg}_\C(\gamma)\ .$$
Hence, $\cS_{g,n}$ is graded,
$$\cS_{g,n} = \bigoplus_{d=0}^{3g-3+n} \cS^d_{g,n}\ .$$
Since the product respects the grading, $\cS_{g,n}$ is a graded
algebra. Of course,
$$q:\cS^d_{g,n} \to H^{2d}(\oM_{g,n},\mathbb{Q})\ . $$

\subsection{The tautological relations $\widetilde{\mathsf{P}}$} 
\label{Ssec:relations}
We define a set $\tP$ consisting of
elements $\cR^d_{g,A}\in \cS^d_{g,n}$ associated to the data
\begin{enumerate}
\item[$\bullet$] $g,n\in \mathbb{Z}_{\geq 0}$ in the stable range $2g-2+n>0$,
\item[$\bullet$] $A=(a_1,\ldots, a_n), \ \ a_i \in\{0,1\}$,
\item[$\bullet$] $d\in \mathbb{Z}_{\geq 0}$ satisfying
$d > \frac{g-1+\sum_{i=1}^n a_i}{3}$.
\end{enumerate}
The elements $\cR^d_{g,A}$ are expressed as sums over
stable graphs of genus $g$ with $n$ legs. 
We prove in Section~\ref{section:cor2} that the conjectured family of relations $\P$ of \cite{Pixton} is implied in cohomology by the family of relations $$q\left(\mathcal{R}^d_{g,A}\right) = 0 \ \ \in H^{2d}(\oM_{g,n},\mathbb{Q})\  $$
for all  $\mathcal{R}_{g,A}^d\in \tP$.
Before writing the formula for $\mathcal{R}^d_{g,A}$, a few definitions are required.

The following two series first arose in the study by Faber and Zagier
of tautological relations on the moduli space $\cM_g$ of nonsingular curves:
\begin{align*}
\B_0(T) &= \sum_{m \geq 0} \frac{(6m)!}{(2m)!(3m)!}(-T)^m
= 1-60T+27720T^2 -\cdots,\\
\B_1(T) &= \sum_{m \geq 0} \frac{1+6m}{1-6m} \frac{(6m)!}{(2m)!(3m)!}(-T)^m
=1 + 84T - 32760T^2 + \cdots.
\end{align*}
These series control the original set $\mathsf{FZ}$ and continue
to play a central role in the set $\tP$.
In the first proof of the Faber-Zagier relations \cite{PPFZ}, 
the above series appeared via differential
equations satisfied by the logarithm of 
$$\Phi(t,x)= \sum_{d=0}^ \infty \prod_{i=1}^d \frac{1}{1-it} \frac{(-1)^d}{d!} 
\frac{x^d}{t^d},$$
see \cite[Section 5]{PPFZ}
and \cite{Ionel}. 
Here we discover a completely different source
for the series $\B_0(T)$ and $\B_1(T)$ 
via the homogeneous calibration of the Frobenius manifold
associated to $A_2$.

Let $f(T)$ be a power series with vanishing constant and linear terms,
$$f(T)\in T^2\mathbb{Q}[[T]]\ .$$
For each $\oM_{g,n}$,  we define
\begin{equation}\label{g33g}
\kappa(f) = \sum_{m \geq 0} \frac1{m!}\ { p_{m*}} \Big(f(\psi_{n+1}) \cdots f(\psi_{n+m})\Big)
 \ \in H^*(\oM_{g,n},\mathbb{Q}),
\end{equation}
where $p_m$ is the forgetful map
$$p_m: \oM_{g,n+m} \to \oM_{g,n}.$$ 
By the vanishing in degrees 0 and 1  of $f$, the sum \eqref{g33g} is finite.

Let $\mathsf{G}_{g,n}$ be the (finite) set of stable graphs of
genus $g$ with $n$ legs (up to isomorphism).
Let $\Gamma \in \mathsf{G}_{g,n}$. For each vertex $v\in \V$,
we introduce an auxiliary variable $\zeta_v$ and impose the
conditions
$$\zeta_v \zeta_{v'}= \zeta_{v'} \zeta_v\ , \ \ \ \zeta_v^2= 1\ .$$
 The variables $\zeta_v$ will be responsible for keeping track of
a local parity condition at each vertex.

The formula for $\mathcal{R}_{g,A}^d$ is a sum over $\mathsf{G}_{g,n}$.
The summand corresponding to $\Gamma \in \mathsf{G}_{g,n}$ is a 
product of 
vertex, leg, and edge factors:
\begin{enumerate}
\item[$\bullet$]
For $v\in \V$, let
$\kappa_v = \kappa\big(T-T \B_0(\zeta_vT)\big)$.
\item[$\bullet$]
For $l \in \L$,
let
 $\B_l = \zeta_{v(l)}^{a_l} \B_{a_l} \! \left(\zeta_{v(l)} \psi_{l}\right)$,  
where $v(l)\in V$ is the vertex to which the leg is assigned.
\item[$\bullet$]
For $e\in \E$, let
\begin{align*}
\Delta_e &= \frac{ \zeta' + \zeta'' - 
\B_0(\zeta' \psi') \zeta''\B_1(\zeta'' \psi'')
-\zeta'\B_1(\zeta' \psi') \B_0(\zeta'' \psi'')}
{\psi'+\psi''}\\
&=(60 \zeta' \zeta''-84) +
\left[32760(\zeta'\psi' + \zeta'' \psi'') - 27720 (\zeta'\psi''+\zeta''\psi')\right] + \cdots,
\end{align*}
where $\zeta',\zeta''$ are the $\zeta$-variables assigned to the vertices adjacent to the edge $e$ and $\psi', \psi''$ are the $\psi$-classes corresponding to the half-edges.
\end{enumerate}
The numerator of $\Delta_e$ is divisible by the denominator due to the identity
$$
\B_0(T) \B_1(-T) + \B_0(-T) \B_1(T) =2.
$$
Obviously $\Delta_e$ is symmetric in the
half-edges.

\begin{definition} \label{Not:relations}
Let $A = (a_1, \dots, a_n) \in \{0,1\}^n$.
We denote by
$\cR_{g,A}^d\in \cS_{g,n}^d$ the degree $d$ component of the strata algebra class 
$$
\sum_{\Gamma\in \mathsf{G}_{g,n}} \frac1{|\Aut(\Gamma)| }
\, 
\frac1{2^{h^1(\Gamma)}}
\;
\left[\Gamma, \; \Bigl[
\prod \kappa_v \prod \B_l 
\prod \Delta_e
\Bigr ]_{\prod_v \zeta_v^{\mathrm{g}(v)-1}}
\right] \ \in \cS_{g,n},
$$
where the products are taken over all vertices, all legs, and all edges of the graph~$\Gamma$.
The subscript $\prod_v \zeta_v^{\mathrm{g}(v)-1}$ indicates
the coefficient of the monomial $\prod_v \zeta_v^{\mathrm{g}(v)-1}$
after the product inside the brackets is expanded.
\end{definition}

\begin{definition}
We denote by $\tP$ the set of classes $\cR^d_{g,A}$ where
$$
d > \frac{g-1 + \sum_{i=1}^n a_i}{3}.
$$
\end{definition}

\begin{theorem} \label{Thm:relations}
Every element $\cR_{g,A}^d\in \widetilde{\mathsf{P}}$ lies in the kernel of 
the homomorphism $$q:\cS_{g,n} \rightarrow H^*(\oM_{g,n},\mathbb{Q})\ .$$ 
\end{theorem}

As a formal consequence of Theorem \ref{Thm:relations}, we will 
establish the originally conjectured
set of relations $\mathsf{P}$. 

\begin{corollary} \label{Pix}
The full set $\mathsf{P}$ of relations conjectured 
in \cite{Pixton} holds in cohomology.
\end{corollary}

Furthermore, we will identify $\cR^d_{g,(a_1, \dots, a_n)}$ as a simple multiple of Witten's class for $r=3$ when $d= \frac{g-1+\sum a_i}{3}$ and a simple multiple of a push-forward of Witten's class under a forgetful map when $d < \frac{g-1 + \sum a_i}{3}$.

\subsection{Cohomological field theories} \label{cft}
We recall here the basic definitions of a cohomological field theory by Kontsevich and Manin~\cite{KonMan}.

Let $V$ be a finite dimensional $\mathbb{Q}$-vector space with 
a non-degenerate symmetric 2-form $\eta$ and a 
distinguished element $\b1 \in V$.
The data $(V,\eta, \b1)$ is the starting point for defining 
a cohomological field theory.
Given a 
basis $\{e_i\}$ of $V$, we write the 
symmetric form as a matrix
$$\eta_{jk}=\eta(e_j,e_k) \ .$$ The inverse matrix is denoted by $\eta^{jk}$ as usual.

A {cohomological field theory} consists of 
a system $\Omega = (\Omega_{g,n})_{2g-2+n > 0}$ of elements 
$$
\Omega_{g,n} \in H^*(\oM_{g,n},\mathbb{Q}) \otimes (V^*)^{\otimes n}.
$$
We view $\Omega_{g,n}$ as associating a cohomology class on $\oM_{g,n}$ 
to elements of $V$ assigned to the $n$ markings.
The CohFT axioms imposed on $\Omega$ are:

\begin{enumerate}
\item[(i)] Each $\Omega_{g,n}$ is $S_n$-invariant, where the action of the symmetric group $S_n$ permutes both the marked points of $\oM_{g,n}$ and the
copies of $V^*$.
\item[(ii)] Denote the basic gluing maps by
$$
q : \oM_{g-1, n+2} \to \oM_{g,n}\ ,
$$
$$
r: \oM_{g_1, n_1+1} \times \oM_{g_2, n_2+1} \to \oM_{g,n}\ .
$$
The pull-backs $q^*(\Omega_{g,n})$ and $r^*(\Omega_{g,n})$ are
 equal to the contractions of $\Omega_{g-1,n+2}$ and 
$\Omega_{g_1, n_1+1} \otimes \Omega_{g_2, n_2+1}$ by the bi-vector 
$$\sum_{j,k} \eta^{jk} e_j \otimes e_k$$ inserted at the two identified points. 
\item[(iii)] Let $v_1, \dots, v_n \in V$ be any vectors and let $p: \oM_{g,n+1} \to \oM_{g,n}$ be the forgetful map. We require 
$$
\Omega_{g,n+1}(v_1 \otimes \cdots \otimes v_n \otimes \b1) = p^*\Omega_{g,n} (v_1 \otimes \cdots \otimes v_n)\ ,
$$
$$\Omega_{0,3}(v_1\otimes v_2 \otimes \b1) = \eta(v_1,v_2)\ .$$
\end{enumerate}

\begin{definition}\label{defcohft}
A system $\Omega= (\Omega_{g,n})_{2g-2+n>0}$ of elements 
$$
\Omega_{g,n} \in H^*(\oM_{g,n},\mathbb{Q}) \otimes (V^*)^{\otimes n}
$$
satisfying properties~(i) and (ii) is a {\em cohomological field theory} or a {\em CohFT}. If (iii) is also satisfied, $\Omega$ is 
a {\em CohFT with unit}.
\end{definition}

A CohFT $\Omega$ yields a {\em quantum product} $\bullet$ on $V$ via 
$$\eta(v_1 \bullet v_2, v_3) = \Omega_{0,3}(v_1 \otimes v_2 \otimes v_3)\ .$$
Associativity of $\bullet$ follows from (ii). The element
$\b1\in V$ is the identity for $\bullet$ by (iii).

A CohFT $\omega$ composed only of degree~0 classes,
$$\omega_{g,n} \in H^0(\oM_{g,n},\mathbb{Q}) \otimes (V^*)^{\otimes n}\ ,$$
 is called a {\em topological field theory}. 
Via property (ii), $\omega_{g,n}(v_1, \dots, v_n)$ 
is determined by considering stable curves with a maximal number of nodes. 
Such a curve is obtained by identifying several rational curves with three marked points. The value of $\omega_{g,n}(v_1 \otimes \cdots \otimes v_n)$ is 
thus uniquely specified by the values of $\omega_{0,3}$ and by the quadratic form~$\eta$. In other words, given $V$ and $\eta$, a topological field theory is uniquely determined by the associated quantum product.

\subsection{Witten's $r$-spin class} \label{wsc}
For every integer $r \geq 2$, there is a beautiful CohFT obtained from 
Witten's $r$-spin class. We review here the basic properties of
the construction. The integer $r$ is fixed once and for all.
  
Let $V$ be an $(r-1)$-dimensional $\mathbb{Q}$-vector space with basis
 $e_0, \dots, e_{r-2}$, bilinear form 
$$
\eta_{ab} = \langle e_a, e_b \rangle =\delta_{a+b,r-2} \, ,
$$
and unit vector $\b1 = e_0$. Witten's $r$-spin theory provides a family of classes
$$
W_{g,n}(a_1, \dots, a_n) \in H^*(\oM_{g,n},\mathbb{Q}).
$$ 
for $a_1, \dots, a_n \in \{0, \dots, r-2 \}$.
These define a CohFT by
$$
\W_{g,n}: V^{\otimes n} \rightarrow H^*(\oM_{g,n},\mathbb{Q}), 
\ \ \ 
\W_{g,n}( e_{a_1} \otimes \cdots \otimes e_{a_n}) =
W_{g,n}(a_1, \dots, a_n)\ .
$$
To emphasize $r$, we will often refer to $V$ as $V_r$.

Witten's class $W_{g,n}(a_1, \dots, a_n)$ has (complex) degree given
by the  formula
\begin{eqnarray}
\label{gred}
\text{deg}_{\C}\ W_{g,n}(a_1, \dots, a_n) & = & 
\D_{g,n}(a_1, \dots, a_n) \\ \nonumber
& = & \frac{(r-2)(g-1) + \sum_{i=1}^n a_i}{r}\ .
\end{eqnarray}
If $\D_{g,n}(a_1, \dots, a_n)$ is not an integer, the corresponding
Witten class vanishes.

In genus 0, the construction was first carried out by Witten \cite{Witten}
using $r$-spin structures ($r^{\rm th}$ roots of the canonical bundle)
and satisfies the following initial conditions:
\begin{equation}\label{fred}
W_{0,3}(a_1,a_2,a_3) = 
\left|
\begin{array}{cl}
1 & \mbox{ if } a_1+a_2+a_3 = r-2,\\
0 & \mbox{ otherwise.}
\end{array}
\right.
\end{equation}
$$
W_{0,4}(1,1,r-2,r-2) = \frac1{r} [\mbox{point}] \ \in H^2(\overline{M}_{0,4},\mathbb{Q})\ .
$$
Uniqueness of Witten's $r$-spin theory in
genus~0 follows easily from the initial conditions \eqref{fred}
and the axioms of a CohFT with unit.

The genus $0$ sector defines a quantum product $\bullet$
on $V$ with unit $e_0$,
$$ \langle e_a\bullet e_b, e_c \rangle = W_{0,3}(a,b,c) \ .$$
The resulting algebra, even after
extension to $\C$, is not semisimple.

The existence of Witten's class in higher genus is both remarkable and
highly non-trivial. 
An algebraic construction was first obtained by Polishchuk and Vaintrob~\cite{PolVai} defining 
$$W_{g,n}(a_1, \dots, a_n)\in A^*(\oM_{g,n}, \mathbb{Q})$$
as an algebraic cycle class. The algebraic approach was 
later simplified by Chiodo~\cite{Chiodo}. 
Analytic constructions have been given by Mochizuki ~\cite{Mochizuki} and
later by Fan, Jarvis, and Ruan \cite{fjr}.
The equivalence between the
above analytic and algebraic constructions was heretofore unknown.

\begin{theorem} \label{Thm:Wittensclass}
For every $r\geq 2$, 
there is a unique CohFT which extends Witten's $r$-spin theory
in genus 0 and has pure dimension \eqref{gred}.
The unique extension takes values in the tautological ring
$$RH^*(\cM_{g,n})\subset H^*(\cM_{g,n},\mathbb{Q}).$$
\end{theorem}

As a consequence of Theorem~\ref{Thm:Wittensclass}, the analytic
and algebraic approaches coincide and yield tautological classes
in cohomology.
Our proof of Theorem~\ref{Thm:Wittensclass} is not valid
for Chow field theories as topological results play an 
essential role.

\subsection{Strategy of proof}
Theorems \ref{Thm:relations} and \ref{Thm:Wittensclass}
are proven together. 
Let $\W_{g,n}$ be any CohFT with unit which
extends Witten's $r$-spin theory in genus 0 and
has pure dimension~\eqref{gred}.

We use a canonical procedure (a shift on the Frobenius manifold) to define a new CohFT $\tW_{g,n}$ satisfying the following four properties:
\begin{enumerate}
\item[(i)] $\tW$ is canonically constructed from $W$ with
 the genus~0 sector of $\tW$ entirely determined by the genus~0 sector of $W$,
\item[(ii)] the quantum product associated to $\tW_{0,3}$ defines a semisimple algebra on $V_r$,
\item[(iii)] the component of $\tW_{g,n}( e_{a_1} \otimes \cdots \otimes e_{a_n})$ in complex degree $\D_{g,n}(a_1, \dots, a_n)$ equals $W_{g,n}(a_1, \dots, a_n)$,
\item[(iv)] the class $\tW_{g,n}( e_{a_1} \otimes \cdots \otimes e_{a_n})$ has no components in degrees higher than $\D_{g,n}(a_1, \dots, a_n)$.
\end{enumerate}
In other words, $\tW$ is constructed from
$\W$ by adding only {\em lower} degree terms.

By the results of Givental and Teleman, $\tW$
is determined via a universal formula in the tautological
ring by the semisimple genus 0 sector.
By property~(iii), we deduce a formula for $\tW$
in the tautological ring depending only upon Witten's $r$-spin
theory in genus 0 and obtain Theorem \ref{Thm:Wittensclass}.

To prove Theorem \ref{Thm:relations}, we write explicitly
Givental's formula for 
the modified CohFT $\tW$ in the $3$-spin case.
The series $\B_0$ and $\B_1$ appear in the associated Frobenius structure. By property~(iv), we obtain vanishings in the tautological ring in degrees
$$
d > \D_{g,n}(a_1, \dots, a_n) =
 \frac{g-1 + \sum_{i=1}^n a_i}{3} \quad \mbox{for } r=3.
$$
The outcome is exactly the relations $\tP$.

As a further outcome of the above investigation, we obtain the
following formula for Witten's 3-spin class.

\begin{theorem} \label{Thm:3Witten}
Let $r=3$. Then, for $g,n \in \Z_{\geq 0}$ in the stable
range, we have
$$W_{g,n}(a_1, \dots, a_n)  = 2^g \, 1728^d \, q\left(\cR^d_{g,(a_1,\ldots, a_n)}\right)
\in H^{2d}(\oM_{g,n},\Q)\ $$
when $d =
 \frac{g-1 + \sum_{i=1}^n a_i}{3}$ is integral (and 
$W_{g,n}(a_1, \dots, a_n)$
is 0 otherwise).\
\end{theorem}

\subsection{Plan of the paper}
In Section~\ref{Sec:FrobA}, we define the shifted Witten class for the
$r$-spin theory. Theorem~\ref{Thm:Wittensclass} is proven as
a consequence of semisimplicity and Teleman's uniqueness result.
A short review of the $R$-matrix action on CohFTs is presented in
Section~\ref{Sec:crashcourse}. 
In Section~\ref{Sec:A2}, 
we compute the $R$-matrix for the $3$-spin case and prove 
Theorems~\ref{Thm:relations} and~\ref{Thm:3Witten}.
The proof of Corollary \ref{Pix} is also given in Section \ref{Sec:A2}.

The study of the $R$-matrix for higher $r$ and the 
 exploration of the associated relations in the tautological ring
will be taken up in \cite{PPZ2}.

\subsection{Acknowledgments} 
We are grateful to A.~Chiodo, P.~Dunin-Barkovsky, C.~Faber,
J.~Gu\'er\'e, F.~Janda, A.~Polishchuk, O.~Randal-Williams, Y. Ruan,
 L.~Spitz, and A.~Vaintrob for useful and detailed discussions. Several 
of the ideas presented here grew out
of discussions at the {\em Geometry and topology of moduli} conference in
Berlin in October 2012 organized by G.~Farkas.
Special thanks to S.~Shadrin for pointing out an error in the first draft of the paper and
to L. Meng (of Peking University) for pointing out a dropped parity factor in the leg term.

R.~P. was partially supported by the Swiss National Science Foundation
grant SNF 200021143274. A.~P. was supported by an NSF Graduate Research Fellowship
and was a guest of the Forschungsinstitut f\"ur Mathematik (FIM) for several visits
to ETH Z\"urich.
D.~Z. was supported by the grant ANR-09-JCJC-0104-01.

\section{$A_{r-1}$ and the shifted Witten class} \label{Sec:FrobA}
\subsection{Potentials}
Frobenius manifolds were introduced and studied in detail in Dubrovin's monograph~\cite{Dubrovin0}. For a concise summary see
\cite[Section 1]{Givental2}.

As for every CohFT, the genus~0 part of Witten's $r$-spin class determines a Frobenius manifold structure
on the underlying vector space $V_r$. 
For Witten's class,
the Frobenius manifold coincides with the 
canonical Frobenius structure on the versal deformation of the $A_{r-1}$ singularity \cite{Dubrovin} up to a coordinate change. We will denote  by $t^0, \dots, t^{r-2}$ the coordinates in the basis $e_0, \dots, e_{r-2}$ of $V_r$. 

The structure of a Frobenius manifold is governed by the
Gromov-Witten potential. The genus~0 Gromov-Witten potential of Witten's 
$r$-spin class 
(without descendants) is:
$$
\F(t^0, \dots, t^{r-2}) = \sum_{n \geq 3} \sum_{a_1, \dots, a_n} \int_{\oM_{0,n}}
W_{0,n}(a_1, \dots, a_n) \frac{t^{a_1} \cdots t^{a_n}}{n!}.
$$
We will refer to $\F$ as the {\em primary genus 0 potential}.

\begin{example}
For $r=3$, the primary genus 0 potential obtained from  Witten's class equals
$$
\F(x,y) = \frac12 x^2 y +\frac1{72} y^4,
$$
where $x = t^0$ and $y = t^1$.

For $r=4$, the potential is 
$$
\F(x,y,z) = \frac12 x^2 z + \frac12 xy^2 + \frac1{16} y^2z^2 + \frac1{960} z^5,
$$
where $x = t^0$, $y = t^1$, and  $z= t^2$.
\end{example}
 The third derivatives of~$\F$ determine an associative algebra structure (the {\em quantum product}) in each tangent space to the Frobenius manifold.
Let $\partial_i$ denote the vector field on $V_r$ associated to 
differentiation by $t^i$. Then,
$$\partial_i \bullet \partial_j = \sum_{k,l}\frac{\partial^3 \F}{\partial t^i \partial t^j
\partial t^k} \eta^{kl} \partial_l\ .$$
The algebra on tangent spaces is semisimple outside the discriminant of $A_{r-1}$. For instance, for $r=3$, the discriminant is $\{ y = 0 \}$.

\begin{definition}\label{ddd}
Let $\tau\in V_{r}$.
We define the {\em shifted Witten class} by
$$
\W_{g,n}^\tau(v_1 \otimes \cdots \otimes v_n) =
\sum_{m \geq 0} \frac1{m!}\ (p_m)_* 
\W_{g,n+m}(v_n\otimes \cdots \otimes v_n \otimes \tau \otimes \cdots \otimes \tau),
$$
where $p_m : \oM_{g,n+m}\to \oM_{g,n}$ is the forgetful map.
\end{definition}

\begin{remark}
We have the following degree bound:
\begin{multline*}
\deg
\Bigl[ (p_m)_* \W_{g,n+m}(e_{a_1} \otimes \cdots \otimes e_{a_n}\otimes
\tau \otimes \cdots \otimes \tau) \Bigr] \\
 \leq \ 
\frac{(g-1)(r-2) + \sum a_i + m(r-2)}{r} -m \\
 = \ \D_{g,n}(a_1, \dots, a_n) - \frac{2m}r\ . 
{ \hspace{+154 pt} }
\end{multline*}
The sum in Definition \ref{ddd} is thus finite for any given 
$g$ and  $a_1, \dots, a_n$. The shifted Witten class is therefore
well-defined. Moreover, the highest degree term of the shifted Witten class 
is equal to the Witten class itself -- all the other terms are of smaller 
degrees. 
\end{remark}

\begin{remark}
Let $\F(t)$ and $\F^\tau(\, \widehat{t}\,)$ be the primary genus~0 potentials of $\W$ and
$\W^\tau$ respectively. By elementary verification, 
$$
\F^\tau(\, \widehat{t} \, ) = \F(\tau + \widehat{t}\,) - (\mbox{\rm terms of degree} <3).
$$
\end{remark}



\begin{proposition}
The shifted Witten class $\W^\tau$ is a CohFT with unit.
\end{proposition}

The proof is a straightforward check, and in any case,
is identical to the proof of Proposition~\ref{Thm:TOmega} given in
 Section \ref{Sec:crashcourse} below.  

\subsection{The Euler field}

A Frobenius manifold is called {\em conformal} if it carries an affine Euler field $E$, a vector field satisfying the following properties:
\begin{enumerate}
\item[(i)]
in flat coordinates $t^i$, the field has the form
$$
E = \sum_i (\alpha_i t^i + \beta_i) \frac\d{\d t^i},
$$ 
\item[(ii)]
the quantum product $\bullet$, the unit $\b1$, and the metric $\eta$ are eigenfunctions of the Lie derivative~$L_E$ with weights $0$, $-1$, and $2-\delta$ respectively. 
\end{enumerate}
The rational number $\delta$ is called the {\em conformal dimension} of the Frobenius manifold. 

For instance, on the Frobenius manifold $A_{r-1}$, an Euler field is given
by 
$$
E = \sum_{a=0}^{r-2} \left(1-\frac{a}{r}\right) t^a \frac{\d}{\d t^a},
$$
$$
\delta = \frac{r-2}{r}.
$$

\begin{remark} We follow here
Givental's conventions for the Euler field. In Teleman's conventions,
the Euler vector field and hence the eigenvalues of $L_E$ have the 
opposite sign.
\end{remark}

Let $\Omega$ be a CohFT and $V$ the corresponding Frobenius manifold. Given an Euler field $E$ on~$V$,
a natural action of $E$ on $\Omega$ is defined
as follows. Let 
$$\deg:H^*(\oM_{g,n},\mathbb{Q}) \to H^*(\oM_{g,n},\mathbb{Q})$$ be the operator which
 acts on $H^{2k}$ by multiplication by~$k$. 
As usual, $\d_i$ is the vector field{\footnote{We
will often use the canonical identification of
$V$ with the tangent space of $0\in V$.}} on $V$ associated to differentiation by
the coordinate $t^i$. 
Then
$$
(E.\Omega)_{g,n}(\d_{i_1}\otimes \dots \otimes \d_{i_n}) = 
$$
$$
\left( \deg + \sum_{l=1}^n \alpha_{i_l} \right) 
\Omega_{g,n}(\d_{i_1} \otimes \dots \otimes \d_{i_n}) 
+ p_* \Omega_{g,n+1} \left(
\d_{i_1} \otimes \dots \otimes \d_{i_n} \otimes  \sum \beta_i \d_i \right),
$$
where $p:\oM_{g,n+1} \to \oM_{g,n}$ is the forgetful map.

\begin{definition}
A CohFT $\Omega$ is {\em homogeneous} if 
$$
(E. \Omega)_{g,n} = [(g-1) \delta +n ]\ \Omega_{g,n}
$$
for all $g$ and $n$.
\end{definition}

Witten's $r$-spin class is easily seen to be homogeneous. Indeed, we have
$$
\D_{g,n}(a_1,\dots,a_n) + \sum_{i=1}^n \left( 1 - \frac{a_i}{r} \right)
=
$$
$$
\frac{(r-2)(g-1) + \sum a_i}{r} + n - \frac{\sum a_i}{r}
=
(g-1) \frac{r-2}{r} + n = (g-1)\delta +n.
$$

The underlying vector space $V_r$ and basis $e_0,\ldots, e_{r-2}$
are the same for the CohFT obtained from the shifted Witten class.
We denote the coordinates on $V_r$ in the basis $e_0,\ldots, e_{r-2}$
for the shifted $r$-spin Witten theory by $\widehat{t}^{\ 0},\ldots, 
\widehat{t}^{\, r-2}$.

\begin{proposition} \label{fvc}
The shifted Witten class is a homogeneous CohFT with Euler field 
$$
E = \sum_{a=0}^{r-2} \left(1-\frac{a}{r}\right) (\tau^a+\widehat{t}^{\, a}) \frac{\d}{\d \widehat{t}^{\, a}},
$$
of conformal dimension $\delta = \frac{r-2}{r}$.
\end{proposition}

\paragraph{Proof.} Assume for simplicity $\tau = u \d_b$ for some fixed $b \in \{ 0, \dots, r-2 \}$. Denote
$$
\W_{g,n+m}
(\d_{a_1} \otimes \dots \otimes \d_{a_n}\otimes \d_b \otimes \cdots \otimes \d_b)
$$
here by just $\W_{g,n+m}$. Then we have
$$
(E.\W^\tau)_{g,n}(\d_{a_1} \otimes \dots \otimes \d_{a_n}) = $$
$$\sum_{m \geq 0} \frac{u^m}{m!}
\left\{
\left[\frac{(r-2)(g-1) + \sum a_i + mb}{r} - m \right]
+ \sum_{i=0}^{n} \left( 1- \frac{a_i}r \right)
\right\}
(p_m)_* \W_{g,n+m}
$$
$$
+ \sum_{m \geq 0} \frac{u^m}{m!}
u \left( 1 - \frac{b}r \right)
(p_{m+1})_* \W_{g,n+m+1}
$$
After simplifying, the above equals
\begin{multline*}
[(g-1)\delta + n] \sum_{m \geq 0} \frac{u^m}{m!} (p_m)_* \W_{g,n+m}\\
- \sum_{m \geq 1} \frac{u^m}{(m-1)!}\left( 1 - \frac{b}r \right) (p_m)_* \W_{g,n+m}
\\
+ \sum_{m \geq 0} \frac{u^{m+1}}{m!}\left( 1 - \frac{b}r \right)  (p_{m+1})_* \W_{g,n+m+1}.
\end{multline*}
The last two sums cancel each other, so we obtain
$$
[(g-1)\delta + n] \sum_{m \geq 0} \frac{u^m}{m!} (p_m)_* \W_{g,n+m}
= [(g-1)\delta + n]\ \W^\tau_{g,n}(\d_{a_1} \otimes \dots \otimes \d_{a_n}) .
$$
The general case is similiar.
\qed

\vspace{10pt}

Since the shifted $r$-spin Witten class is 
a homogeneous semisimple CohFT,
 we can apply the following theorem by C.~Teleman \cite[Theorem 1]{Teleman}.

\begin{theorem}[Teleman] \label{tele}
Let $\Omega_{0,n}$ be a genus 0 homogeneous semisimple Coh\-FT with unit.
The following results hold:
\begin{enumerate}
\item[(i)] There exists a unique homogeneous CohFT with unit $\Omega_{g,n}$ extending
$\Omega_{0,n}$ to higher genus.
\item[(ii)] The extended CohFT $\Omega_{g,n}$ is obtained
by an $R$-matrix action on the topological (degree 0) sector of 
$\Omega_{0,n}$ determined by $\Omega_{0,3}$.
\item[(iii)] The $R$-matrix is uniquely specified by  $\Omega_{0,3}$ and 
the Euler field.
\end{enumerate}
\end{theorem}


\noindent The unit-preserving $R$-matrix action in part (ii) of Theorem \ref{tele} will be reviewed in Section~\ref{Sec:crashcourse}, see
Definition~\ref{Def:RactionWithUnit}. 

We will compute the $R$-matrix for the 3-spin Witten class in Section~\ref{Sec:A2}. Since the shifted $3$-spin Witten class (considered for all genera)
is a homogeneous CohFT with unit,
 the expressions obtained by the $R$-matrix action 
coincide with the shifted $3$-spin Witten class. 
In particular, if we split the expression of the $R$-matrix action into pure degree parts, the parts of degree $$d > \D_{g,n}(a_1, \dots, a_n)$$ vanish while the
part of degree~$\D_{g,n}(a_1, \dots, a_n)$ coincides with Witten's class, which proves Theorem~\ref{Thm:Wittensclass}.

\section{The $R$-matrix action} \label{Sec:crashcourse}

We present here a succinct but self-contained review of the $R$-matrix action on CohFTs. The action was first defined on Gromov-Witten potentials by Givental~\cite{Givental}. Its lifting to CohFTs was independently discovered by several authors:  the papers by Teleman~\cite{Teleman} and Shadrin~\cite{Shadrin} give 
an abbreviated 
treatment of the subject and refer to unpublished notes by Kazarian and by Katzarkov, Kontsevich, and Pantev.

%
%

\subsection{The $R$-matrix action on CohFTs}
\label{Ssec:Rnounit}
Let $V$ be a vector space with basis $\{e_i\}$ and 
a symmetric bilinear form $\eta$.
Consider the group of $\End(V)$-valued power series 
\begin{equation}\label{gddx}
R(z) = 1 + R_1 z + R_2 z^2 + \cdots
\end{equation}
satisfying the {\em symplectic condition}, 
$$R(z)R^*(-z) = 1\ ,$$ 
where $R^*$ is the adjoint with respect to $\eta$.

\begin{remark}
Let $R^k_j$ be the matrix form of an endomorphism $R$ in
the given basis,
$$ R\left({t^j e_j}\right) = \sum_{j,k} R^k_j t^j \, e_k\ .$$
The symplectic condition in coordinates is
$$
\sum_{l,s,k} 
R_l^j(z) \eta^{ls} R_s^k(-z) \eta_{ku} = \delta^j_u.
$$
After multiplying by $\eta^{-1}$ on the right,
we obtain an equivalent condition in bi-vector form
$$
\sum_{l,s} R_l^j(z) \eta^{l s} R_s^k(-z)= \eta^{jk}\ .
$$
We conclude that the expression
$$
\frac{\eta^{jk}-\sum_{l,s} R_l^j(z) \eta^{l s} R_s^k(w)}{z+w}
$$
is a well-defined {\em power series} in $z$ and~$w$. 

Associated to $R(z)$ is the power series $R^{-1}(z)= \frac{1}{R(z)}$ which
also satisfies the symplectic condition.{\footnote{By the symplectic condition, we have
$R^{-1}(z) = R^*(-z)$.}} It follows that
\begin{equation} \label{Eq:zw}
\sum_{j,k} \frac{\eta^{jk}-\sum_{l,s} (R^{-1})_l^j(z) \eta^{l s} (R^{-1})_s^k(w)}{z+w} e_j \otimes e_k
\in V^{\otimes 2}[[z,w]].
\end{equation}
\end{remark}
We will denote the $V^{\otimes 2}$-valued power series \eqref{Eq:zw} by
$$\frac{\eta^{-1} - R^{-1}(z) \eta^{-1} R^{-1}(w)^\xxx}{z+w} $$
for short (where the superscript $\xxx$ denotes matrix transpose).

Let $\Omega = (\Omega_{g,n})_{2g-2+n>0}$ be a CohFT on~$V$,
 and let $R$ be an element of the group \eqref{gddx}. 
The CohFT $R \Omega$ is defined as follows.

\begin{definition} \label{rrr}
Let $\mathsf{G}_{g,n}$ be the finite set of stable graphs{\footnote{See
Sections \ref{stgr}-\ref{Ssec:relations}.}}
 of genus $g$
with $n$ legs.
For each $\Gamma\in \mathsf{G}_{g,n}$, define a  contribution
$$\text{Cont}_{\Gamma} \in  H^*(\oM_{g,n},\mathbb{Q}) \otimes (V^*)^{\otimes n}$$
by the following construction:
\begin{enumerate}
\item[(i)]
place $\Omega_{\mathsf{g}(v),\mathsf{n}(v)}$ at each vertex $v$ of~$\Gamma$,
\item[(ii)]
place $R^{-1}(\psi_l)$ at every leg~$l$ of~$\Gamma$,
\item[(iii)] at every edge $e$ of $\Gamma$, place
$$
\frac{\eta^{-1}-R^{-1}(\psi'_e) \eta^{-1} R^{-1}(\psi''_e)^{\xxx}}{\psi_e' + \psi_e''}\ .$$
\end{enumerate}
Define $(R\Omega)_{g,n}$ to be the sum of contributions of all stable graphs,
$$(R\Omega)_{g,n} = \sum_{\Gamma \in \mathsf{G}_{g,n}} \frac{1}{|\Aut(\Gamma)|} \ 
\text{Cont}_\Gamma\ .$$ 
\end{definition}

We use the inverse of the $R$-matrix in all of our formulas in
Definition \ref{rrr}. 
There are two reasons for the
seemingly peculiar choice.
First, the result will be a  left group action rather than a right group action on CohFTs. Second,
the same convention is used by Givental and Teleman in their papers.

A few remarks about Definition \ref{rrr} are needed for clarification.
By the symmetry property of CohFTs, the placement of $\Omega_{\mathsf{g}(v),
\mathsf{n}(v)}$ does not depend upon an ordering of the
half-edges at $v$. 
At a leg $l$ attached to a vertex~$v$, we have
$$R^{-1}(\psi_l) \in H^*(\oM_{\mathsf{g}(v),\mathsf{n}(v)},\mathbb{Q}) 
\otimes \End(V).$$ 
The first factor acts on the cohomology of the moduli space 
$\oM_{\mathsf{g}(v), \mathsf{n}(v)}$ by multiplication. 
The endomorphism factor acts
on the vectors  
which are 
``fed'' to $\Omega_{\mathsf{g}(v),\mathsf{n}(v)}$
at the legs. 

For an edge $e$ attached to vertices $v'$ and $v''$ (possibly the same vertex), denote by $\oM_{g',n'}$ and $\oM_{g'',n''}$ the 
corresponding moduli spaces. The insertion on~$e$ is an element of 
$$H^*(\oM_{g',n'},\mathbb{Q}) \otimes H^*(\oM_{g'',n''},\mathbb{Q}) 
\otimes V^{\otimes 2}\ $$
obtained by substituting $z = \psi_e'$ and $w = \psi_e''$ in (\ref{Eq:zw}). 
Once again, the cohomology factors act on the corresponding cohomology spaces by multiplication. The bivector part is used to contract the two covectors sitting on the half-edges $e'$ and $e''$ in the corresponding CohFT elements
at $v'$ and $v´´$. In the expression 
$R^{-1}(\psi_e') \eta^{-1} R^{-1}(\psi_e'')^{\xxx}$, the bivector $\eta^{-1}$ sits in the middle of the edge, while the action of $R^{-1}$ is directed from the middle of the edge towards the vertices. 

The similarity of Definition \ref{rrr}  with the form of the relations~$\cR^d_{g,A}$ was the starting point of our paper.

\begin{proposition}
If $\Omega$ is a CohFT, 
the system $(R\Omega)_{g,n}$ is a CohFT.
\end{proposition}

\paragraph{Proof.} 
The symmetry of $(R\Omega)$ follows directly from the
symmetry of $\Omega$ and the definition of the $R$-matrix action. Hence, we need only
establish the pull-back property (ii) of Definition \ref{defcohft}.

Let $\Phi \in \mathsf{G}_{g,n}$ be a stable graph with a
single edge $e$.
In order to compute the pull-back of tautological classes
under $$\xi_\Phi: \oM_\Phi \rightarrow \oM_{g,n}\ ,$$
according to the rule given in Section~\ref{Ssec:introdualgraphs},
 we must enumerate all stable graphs $\Gamma$ 
with a distinguished edge $e \in E(\Gamma)$ such that 
contracting all other edges 
yields the graph~
$\Gamma_2 = \Phi$.

If $E_1 = E$, we have $\Gamma_1 = \Gamma$. 
The contribution of $\Gamma$ to the  pull-back is obtained from the contribution of $\Gamma_1$ to $R\Omega$ after a multiplication by 
$$-(\psi_e' + \psi_e'')\ .$$ 
In other words, the contribution to the pull-back is 
obtained by placing the class 
$$R^{-1}(\psi'_e) \eta^{-1} R^{-1}(\psi''_e)^\xxx - \eta^{-1}$$ 
on the edge~$e$ with the usual insertions on all other edges and legs.

If $E_1 = E \setminus \{ e \}$ then $\Gamma_1$ is obtained from $\Gamma$ by contracting~$e$. According to the CohFT rules for $\Omega$, the contribution of $\Gamma$ to the pull-back is obtained by placing $\eta^{-1}$ on the edge~$e$ with the standard classes on all other edges.

After summing, the total contribution of $\Gamma$ to the pull-back is equivalent to placing 
$$
R^{-1}(\psi'_e) \, \eta^{-1} \,  R^{-1}(\psi''_e)^\xxx.
$$ 
on the edge $e$ of $\Gamma$ --
precisely what we obtain by the CohFT rules applied to $(R\Omega)_{g,n}$.
\qed

\vspace{+8pt}
\begin{proposition}
The $R$-matrix action on CohFTs  is a left group action. 
\end{proposition}

\paragraph{\em Proof.} We must prove the action of $R_a(z)$ followed by the
 action of $R_b(z)$ is equal to the action of $R_b(z) R_a(z)$.

When we apply $R_a(z)$, we sum over all stable graphs of type $(g,n)$. Let us color the edges of these stable graphs in red. When we then apply $R_b(z)$,
 we sum over all stable graphs of type $(g,n)$, but now we replace each vertex of the stable graph by a small red graph. Let us color the edges of the large graph in blue.
 The result of the consecutive actions of $R_a(z)$ and $R_b(z)$ will be a sum over all stable graphs of type $(g,n)$ whose edges are colored in red and blue. 

On every leg~$l$ of the stable graph, we place first $R^{-1}_a(\psi_l)$ closer to the vertex and then $R^{-1}_b(\psi_l)$ at the end of the leg. The final outcome is 
\begin{equation}\label{nndd}
R^{-1}_a(\psi_l) R^{-1}_b(\psi_l) = (R_b R_a)^{-1}(\psi_l)\ .
\end{equation}
The result \eqref{nndd}
 is also what we place on a leg when we compute the action of the product $R_b R_a$.

Consider next an edge~$e$ of the stable graph. 
 We will use the abbreviations  $R' = R(\psi'_e)$ and $R'' = R(\psi''_e)$. 
On a red edge, we have placed
$$
\mbox{\em red edge:} \qquad
\frac{\eta^{-1}- (R'_a)^{-1} \eta^{-1} ((R'_a)^{-1})^\xxx}{\psi_e' + \psi_e''}
$$
via the first action. On a blue edge, on the other hand, we see  $R'_a$ and $R''_a$ on the ends of the edge and  
$$
\frac{\eta^{-1}-(R'_b)^{-1} \eta^{-1} ((R'_b)^{-1})^\xxx}{\psi_e' + \psi_e''}.
$$
in the middle of the edge. The final outcome, after unwinding the definitions, is
$$
\mbox{\em blue edge:} \qquad
\frac{(R'_a)^{-1} \eta^{-1} ((R'_a)^{-1})^\xxx - (R'_a)^{-1}(R'_b)^{-1} \eta^{-1} 
((R'_b)^{-1})^\xxx
((R'_a)^{-1})^\xxx}{\psi_e' + \psi_e''}.
$$
Since we are summing over all possible colorings, every edge in the stable graph will appear once in red and once in blue. The total contribution will be the sum of the contributions of the two colors, 
$$
\mbox{\em red } + \mbox{ \em blue:} \qquad
\frac{\eta^{-1} - (R'_bR'_a)^{-1} \eta^{-1} ((R'_bR'_a)^{-1})^\xxx}{\psi_e' + \psi_e''}.
$$
The result is exactly what is placed on an edge when we compute the action of the product $R_b R_a$. \qed

\subsection{Action by translations} \label{Ssec:translations}

Let $\Omega$ be a CohFT based on the vector space~$V$, and let $$
T(z) = T_2 z^2 + T_3z^3 + \cdots
$$ 
be a $V$-valued power series with vanishing  coefficients in degrees 0 and 1.

\begin{definition}
The {\em translation of $\Omega$ by~$T$} is the CohFT $T\Omega$ defined by
$$
(T\Omega)_{g,n}(v_1 \otimes \cdots \otimes v_n)
$$
$$
= 
\sum_{m \geq 0} \frac1{m!}\ (p_m)_* \Omega_{g,n+m} 
\Bigl(
v_1 \otimes \cdots \otimes v_n \otimes T(\psi_{n+1}) \otimes \cdots \otimes T(\psi_{n+m})
\Bigr),
$$
where $p_m : \oM_{g,n+m}\to \oM_{g,n}$ is the forgetful map.
\end{definition}

The use of  $T(\psi_i)$ as an argument in a CohFT is an abuse of notation.
The result should be understood as
$$
\Omega_{g,n}( \cdots T(\psi_i) \cdots) = \sum_{k \geq 2} \psi_i^k 
\Omega_{g,n}( \cdots T_k \cdots) .
$$

\begin{remark}
The action by translations is very close to the shift of
Definition \ref{ddd}. However, unlike shifts, the translation action is always well-defined for degree reasons: the degree of the $m^{\rm th}$ summand of the definition is at least~$m$, so the sum is actually finite for any given~$g,n$.
\end{remark}

The action by translations can be described in terms of stable graphs. 
It is a summation over stable graphs with a single vertex 
and $n+m$ legs for $m \geq 0$. The first $n$ legs carry the vectors 
$v_1, \dots, v_n$, and the last $m$ legs carry the series $T(\psi_i)$. 
The latter legs are then suppressed by a forgetful map. 
We will call the first $n$ legs {\em main legs} and the last $m$ legs {\em $\kappa$-legs}, since the push-forward of powers of $\psi$-classes gives rise to $\kappa$-classes.

\begin{proposition} \label{Thm:TOmega}
If $\Omega$ is a CohFT, the system 
$(T\Omega)_{g,n}$ is a CohFT.
\end{proposition}

\paragraph{Proof.} Let $\Phi$ be a stable graph with a
single edge, and let $\oM_\Phi$ the corresponding moduli space. 
We examine the pull-back of $T\Omega$ to $\oM_\Phi$. 
If $\Phi$ has a single vertex, then the $\kappa$-legs in the definition of $T\Omega$ just stay on this vertex. If $\Phi$ has two vertices, then the $\kappa$-legs are distributed among the two vertices. The automorphism coefficients match: there are 
$\binom{m}{m_1,m_2}$ ways to distribute $m$ $\kappa$-legs between two vertices, which leads to a coefficient 
$$
\frac{1}{m!} \binom{m}{m_1,m_2} = \frac1{m_1!} \, \frac1{m_2!}\ .
$$
Thus, $T\Omega$ satisfies the axioms of a CohFT. \qed

\begin{proposition} \label{bnnb}
Translations form an abelian group action on CohFTs.
\end{proposition}

\paragraph{Proof.} The definition of $(T_a + T_b)\Omega$ contains the following sum (where we have suppressed the $\psi$-classes in the notation):
$$
\sum_{m\geq 0} \frac{(T_a+T_b)^m}{m!}  = 
\sum_{m \geq 0} \sum_{m_a+m_b=m} \frac{T_a^{m_a} T_b^{m_b}}{m_a! \, m_b!}  =
\sum_{m_a \geq 0} \frac{T_a^{m_a}}{m_a!} \sum_{m_b \geq 0} \frac{T_b^{m_b}}{m_b!},
$$
which is the definition of the successive actions of $T_b$ and $T_a$. \qed

\begin{proposition} \label{Thm:TRequalsRT} 
Let $R(z) \in {\rm id} + z\End(V)[[z]]$ be an $\End(V)$-valued power series satisfying the symplectic condition. Let $T_a, T_b \in z^2 V[[z]]$ be two $V$-valued power series related to $R$ by 
$$T_a(z) = R(z) T_b(z)\ .$$ Then, for every CohFT $\Omega$, we have
$$
T_a R \Omega = R T_b \Omega\ .
$$
\end{proposition}

\begin{remark} \label{Rem:affinegroup}
The equality of the proposition can be written more concisely as
$$
RT\Omega = (RT)R\Omega
$$
or
$$
R T R^{-1} \Omega = (RT) \Omega.
$$
Thus the actions of $R$ and $T$ can be combined into an action of an affine group.
\end{remark}

\paragraph{Proof of Proposition~\ref{Thm:TRequalsRT}.} Both CohFTs $T_a R \Omega$ and $RT_b \Omega$ can be expressed as sums over stable graphs. We will
match the sums.

Consider first $T_a R \Omega$. By definition, 
we start with $n$ main legs marked by $v_1, \dots, v_n$ and $m$ $\kappa$-legs 
marked by $T_a(\psi_{n+1})$, \dots, $T_a(\psi_{n+m})$. 
We attach these legs to all possible stable graphs of genus $g$ with $n+m$ legs. Their vertices are marked with $\Omega$, their legs with $R^{-1}(\psi_i)$, and their edges with 
$$
\frac{\eta^{-1} - R^{-1}(\psi') \eta^{-1} R^{-1}(\psi'')^\xxx}{\psi' + \psi''}.
$$
The outcome is a sum over stable graphs of genus $g$
with $n+m$ legs whose $n$ main legs are marked with $R^{-1}(\psi_i)(v_i)$ and $m$ $\kappa$-legs with $$(R^{-1} T_a)(\psi_i) = T_b(\psi_i)\ .$$

Consider next $R T_b \Omega$. 
 We start with a sum over all stable graphs of genus $g$
 with $n$ legs. Their vertices are marked with $\Omega$, their legs with $R^{-1}(\psi_i)$, and their edges with 
$$
\frac{\eta^{-1} - R^{-1}(\psi') \eta^{-1} R^{-1}(\psi'')^\xxx}{\psi' + \psi''}.
$$
Now we add to every vertex of such a graph an arbitrary number of $\kappa$-legs marked with $T_b(\psi_i)$. The sum only runs over the graphs which remain stable when we remove the $\kappa$-legs. However, the summation can be extended to all stable graphs with $n+m$ legs.
Indeed, if a stable graph has a genus~0 vertex~$v$ with $m$ $\kappa$-legs and less than 3 other half-edges, then the dimension of the moduli space assigned to~$v$ is less than $m$, while the degree of the class sitting on this moduli space is at least~$2m$. Thus the contribution of the graph vanishes. In conclusion, 
we obtain exactly the same sum as in the first case. 

The sums here are infinite (as $m$ is unbounded), but only a finite number 
of terms are nonzero. The same issue arose in the definition of the translation action.
\qed

\subsection{CohFTs with unit} \label{Ssec:withunit}

Let $\Omega$ be 
a CohFT with unit~$\b1 \in V$ satisfying
$$
\Omega_{g,n+1}(v_1 \otimes \cdots \otimes v_n \otimes \b1) =
p^* \Omega_{g,n}(v_1 \otimes \cdots \otimes v_n), \
$$
$$ \Omega_{0,3}(v_1 \otimes v_2 \otimes \b1) = \eta(v_1,v_2) \ ,$$
where $p:\oM_{g,n+1} \to \oM_{g,n}$ is the forgetful map.
The $R$-matrix action and the translation action defined in the previous sections do not preserve the property of being a CohFT with unit. 
However, we will explain here how the
 two actions can be combined in a unique way so as to preserve 
the unit. 

We recall a well-known geometric result
 which  we will implicitly use in the computations.
\begin{lemma}
Consider the following commutative square of forgetful maps:
\begin{center}
\begin{picture}(120,100)
\put(0,0){$\oM_{g,n+k}$}
\put(100,0){$\oM_{g,n}$}
\put(0,80){$\oM_{g,n+k+m}$}
\put(100,80){$\oM_{g,n+m}$}
\put(10,70){\vector(0,-1){55}}
\put(110,70){\vector(0,-1){55}}
\put(43,5){\vector(1,0){50}}
\put(58,85){\vector(1,0){35}}
\put(15,40){$P_m$}
\put(90,42){$p_m$}
\put(65,70){$P_k$}
\put(65,10){$p_k$}
\end{picture}
\end{center}
The relation  
$(p_k)^* (p_{m})_* = (P_{m})_* (P_k)^*$ holds in cohomology.
\end{lemma}
\paragraph{Proof.}
Let $X$ be the fiber product of $p_m$ and $p_k$, with maps 
$$a:X\to\oM_{g,n+m}\ , \ \ b:X\to\oM_{g,n+k}\ ,  \ \ f:\oM_{g,n+k+m}\to X\ .$$
 Then $(p_k)^*(p_m)_* = b_*a^*$ is immediate, and also
\[
(P_m)_*(P_k)^* = (b_* f_*)(f^* a^*) = b_* (f_* f^*) a^* = b_* a^*
\]
by birationality of $f$. \qed
\vspace{10pt}

The definition of the translation action involves push-forwards and the axioms of a CohFT with unit involve a pull-back. By the above Lemma,
 we will not have to worry about whether the pull-back is taken before or after the push-forward.

\begin{proposition} 
Let $\Omega$ be a CohFT with unit $\b1\in V$.
Let $R(z)$ be an $R$-matrix satisfying the symplectic condition, and let
 $$T_a(z) = z \cdot [R(\b1) - \b1](z),\ \ T_b(z) = z \cdot [\b1 - R^{-1}(\b1)](z)
$$
be two elements of $z^2V[[z]]$.
Then,
$$T_aR\Omega = RT_b\Omega
$$
is also a CohFT with unit $\b1\in V$.
\end{proposition}

\paragraph{Proof.}  By Proposition~\ref{Thm:TRequalsRT}, $T_aR\Omega$
and $RT_b\Omega$ are equal CohFTs. Hence, only the unit property
must be verified.

The CohFT $RT_b\Omega_{g,n}(v_1 \otimes \cdots \otimes v_n)$ is expressed as
 a sum over stable graphs of genus $g$ with $n$ main legs and any number $m \geq 0$ of $\kappa$-legs (see the proof of Proposition~\ref{Thm:TRequalsRT}).
Their vertices are marked with $\Omega$, their main legs with $R^{-1}(\psi_i)$, 
their $\kappa$-legs with $T_b(\psi_i)$, and their edges with 
$$
\frac{\eta^{-1} - R^{-1}(\psi') \eta^{-1} R^{-1}(\psi'')^\xxx}{\psi' + \psi''}.
$$
The expression is an infinite sum with only a finite number of nonzero terms. 
To calculate $p^*RT_b\Omega_{g,n}(v_1 \otimes \cdots \otimes v_n)$, we will study the pull-back under $p^*$ of the contribution of 
every stable graph $\Gamma$ in the
 sum.

Let us call the new  leg marked $n+1$
appearing on $\oM_{g,{n+1}}$ after the pull-back the {\em special leg}.
The pull-back of a stable graph $\Gamma$ is given
by the stable graphs obtained by attaching the special leg to one of the vertices of~$\Gamma$. 

The pull-backs of the stable graph contributions involve also
the pull-backs of the cotangent line classes.
The relation between the pulled-back $\psi$-classes $p^* \psi_i$ from
$\oM_{g,n+m}$ in terms of the new classes $\psi_i$ on $\oM_{g,n+m+1}$ is given by the well-known formula
\begin{equation} \label{Eq:pullbackpsi}
p^*(\psi_h^d) = \psi_h^d - \Delta_{h,n+1}\, p^*(\psi_h^{d-1}).
\end{equation}
Here $h$ is a half-edge of~$\Gamma$ and
$\Delta_{h,n+1}$ is the divisor{\footnote{The divisor is empty
unless $h$ and the special leg are on the same vertex.}} in the moduli space at the vertex 
carrying $h$ corresponding to 
curves with a genus 0 component carrying only the markings $h$ and $n+1$.

The pull-back of the contribution of
$\Gamma$ is given by a sum of two kinds of terms. 
A term of the first kind is obtained by attaching the special leg to one of the vertices of $\Gamma$ and placing the class $\b1$ on it. This 
 happens if we choose the first term on the
right side of (\ref{Eq:pullbackpsi}) in the pull-back of each $\psi$-class.
We use here also the original unit property of $\Omega$.

A term of the second kind is obtained by choosing a half-edge of $\Gamma$, placing a new vertex on it, and adding the special leg maked with $\b1$ to this vertex. The power of the $\psi$-class that was written on this half-edge is then reduced by~1. This happens if we choose the second term of (\ref{Eq:pullbackpsi}) in the pull-back of the $\psi$-class corresponding to the chosen half-edge. 
A term of the second kind occurs with a minus sign.

Let us look more closely at the terms of the second kind.
Suppose we have placed the new vertex on the $i^{\rm th}$ main leg. Then there are two legs attached to this vertex. 
\begin{center}
\ 
\begin{picture}(0,0)%
\includegraphics{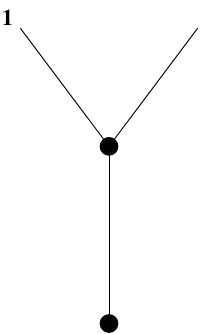}%
\end{picture}%
%
%
\setlength{\unitlength}{3947sp}%
\begingroup\makeatletter\ifx\SetFigFont\undefined%
\gdef\SetFigFont#1#2#3#4#5{%
  \reset@font\fontsize{#1}{#2pt}%
  \fontfamily{#3}\fontseries{#4}\fontshape{#5}%
  \selectfont}%
\fi\endgroup%
\begin{picture}(1647,2666)(2436,-2783)
\put(3376,-2503){\makebox(0,0)[lb]{\smash{{\SetFigFont{12}{14.4}{\rmdefault}{\mddefault}{\updefault}{$\psi' = p^* \psi_i$}%
}}}}
\put(3393,-1661){\makebox(0,0)[lb]{\smash{{\SetFigFont{12}{14.4}{\rmdefault}{\mddefault}{\updefault}{$\psi''$}%
}}}}
\put(3543,-1128){\makebox(0,0)[lb]{\smash{{\SetFigFont{12}{14.4}{\rmdefault}{\mddefault}{\updefault}{$\psi_i$}%
}}}}
\put(2659,-1094){\makebox(0,0)[lb]{\smash{{\SetFigFont{12}{14.4}{\rmdefault}{\mddefault}{\updefault}{$\psi_{n+1}$}%
}}}}
\put(4068,-311){\makebox(0,0)[lb]{\smash{{\SetFigFont{12}{14.4}{\rmdefault}{\mddefault}{\updefault}{$v_i$}%
}}}}
\end{picture}%

\end{center}

First, the $i^{\rm th}$ main leg carrying a $v_i$ is attached. 
We can replace the $v_i$ by $R^{-1}(v_i)(\psi_i)$ since anyway $\psi_i=0$ on the moduli space $\oM_{0,3}$ corresponding to our vertex. Second, the special leg 
carrying $\b1$ is attached. We can similarly replace $\b1$ by 
$R^{-1}(\b1)(\psi_{n+1})$. Finally, there is the edge connecting our vertex to the rest of the graph with a 
\begin{equation}\label{fnnv}
\frac{[v_i-R^{-1}(v_i)](p^*\psi_i)}{p^*\psi_i}
\end{equation}
placed{\footnote{The minus sign in
the second term of \eqref{Eq:pullbackpsi} is included in \eqref{fnnv}.}} on it. We can replace the edge insertion with the standard one
$$
\frac{\eta^{-1} - R^{-1}(\psi') \eta^{-1} R^{-1}(\psi'')^\xxx}{\psi' + \psi''}
$$
since we have $\psi''=0$ and 
$$\Omega_{0,3}(v_i \otimes v_i'' \otimes \b1) = \eta(v_i, v_i'').$$
We have therefore obtained a stable graph with $n+1$ legs marked precisely 
as in the definition of 
$(RT_b\Omega)_{g,n+1}(v_1 \otimes \cdots \otimes v_n \otimes \b1)$.

Next, suppose we have placed the new vertex on a half-edge. We
group the terms obtained from  the two half-edges of
a single edge of $\Gamma$ together.

\begin{center}
\ 
\begin{picture}(0,0)%
\includegraphics{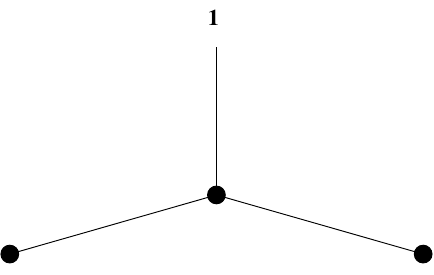}%
\end{picture}%
%
%
\setlength{\unitlength}{3947sp}%
\begingroup\makeatletter\ifx\SetFigFont\undefined%
\gdef\SetFigFont#1#2#3#4#5{%
  \reset@font\fontsize{#1}{#2pt}%
  \fontfamily{#3}\fontseries{#4}\fontshape{#5}%
  \selectfont}%
\fi\endgroup%
\begin{picture}(3465,2198)(1812,-2399)
\put(3609,-1353){\makebox(0,0)[lb]{\smash{{\SetFigFont{12}{14.4}{\rmdefault}{\bfdefault}{\updefault}{$\psi_{n+1}=0$}%
}}}}
\put(1993,-2019){\makebox(0,0)[lb]{\smash{{\SetFigFont{12}{14.4}{\rmdefault}{\bfdefault}{\updefault}{$\psi'$}%
}}}}
\put(4593,-2328){\makebox(0,0)[lb]{\smash{{\SetFigFont{12}{14.4}{\rmdefault}{\bfdefault}{\updefault}{$\psi''$}%
}}}}
\put(2618,-1786){\makebox(0,0)[lb]{\smash{{\SetFigFont{12}{14.4}{\rmdefault}{\bfdefault}{\updefault}{$\psi'''=0$}%
}}}}
\put(3576,-2136){\makebox(0,0)[lb]{\smash{{\SetFigFont{12}{14.4}{\rmdefault}{\bfdefault}{\updefault}{$\psi''''=0$}%
}}}}
\end{picture}%

\end{center}

\noindent The standard edge insertion for $\Gamma$ is 
$$
\frac{\eta^{-1} - R^{-1}(\psi') \eta^{-1} R^{-1}(\psi'')^\xxx}{\psi' + \psi''}\ .
$$
If we place the new vertex at first half-edge, we obtain
$$
- \frac1{\psi'} 
\left[
\frac{\eta^{-1} - R^{-1}(\psi') \eta^{-1} R^{-1}(\psi'')^\xxx}{\psi' + \psi''}
- \frac{\eta^{-1} - \eta^{-1} R^{-1}(\psi'')^\xxx}{\psi''}
\right]\ .
$$
Here, we have subtracted the $\psi'$-free term from the edge insertion and divided the result by $\psi'$. 
The minus sign in front is the sign of the second term of~\eqref{Eq:pullbackpsi}.
Similarly, if we place the new vertex on the second half-edge, we obtain
$$
- \frac1{\psi''} 
\left[
\frac{\eta^{-1} - R^{-1}(\psi') \eta^{-1} R^{-1}(\psi'')^\xxx}{\psi' + \psi''}
- \frac{\eta^{-1} - R^{-1}(\psi')\eta^{-1}}{\psi'}
\right]\ .
$$
Adding the two contributions yields
$$
\frac{\eta^{-1} - R^{-1}(\psi')\eta^{-1} - \eta^{-1} R^{-1}(\psi'')^\xxx
+ R^{-1}(\psi') \eta^{-1} R^{-1}(\psi'')}{\psi' \psi''}^\xxx =
$$
\begin{equation}\label{hjjg}
\frac{\eta^{-1} - R^{-1}(\psi')\eta^{-1}}{\psi'} 
\;\; \eta \;\;
\frac{\eta^{-1} - \eta^{-1}R^{-1}(\psi')^\xxx}{\psi''}\ .
\end{equation}
The result \eqref{hjjg}
 is precisely the product of the standard edge insertions for the two new edges, considering the $\psi$-classes at the new vertex vanish and 
$$\Omega_{0,3}(v' \otimes v'' \otimes \b1) = \eta(v', v'')\ .$$
 Finally, as before, we replace the $\b1$ on the special leg by $R^{-1}(1)(\psi_{n+1})$ without consequence since $\psi_{n+1}=0$ on our vertex. 
Once again, we have obtained a stable graph with $n+1$ legs marked precisely as in the definition of $(RT_b\Omega)_{g,n+1}(v_1 \otimes \cdots \otimes v_n \otimes \b1)$.

The final case to consider is when we place the new vertex on a $\kappa$-leg. 
The edge joining the new vertex to the rest of the graph will then be marked by 
$$
-\frac{T_b(\psi_i)}{\psi_i} = [R^{-1}(\b1) - \b1](\psi_i).
$$
We immediately take the push-forward of our class under the partial forgetful 
map that forgets just the single $\kappa$-leg we are considering. 
We will obtain a graph on which the special leg carries the marking 
$$
[R^{-1}(\b1) - \b1](\psi_{n+1}).
$$
Exactly the same graph also appears among what we called the terms of the first kind --- when we attach the special leg without creating any new vertices. 
There the special leg carried the marking $\b1$. After 
adding the two contributions together,  we obtain
 $R^{-1}(\b1)(\psi)$, which is the standard leg insertion.

We have shown $p^*(RT_b \Omega)_{g,n}(v_1 \otimes \cdots \otimes v_n)$ 
is given by precisely the same sum over stable graph contributions as
$(RT_b \Omega)_{g,n+1}(v_1 \otimes \cdots \otimes v_n \otimes \b1)$.
 Therefore the two are equal. \qed

\begin{definition} \label{Def:RactionWithUnit}
Let $\Omega$ be a CohFT with unit $\b1\in V$.
Let $R(z)$ be an $R$-matrix satisfying the symplectic condition, and let
 $$ T(z) = z \cdot \b1 - zR^{-1}(z)(\b1) \in z^2V[[z]]. \ $$
The {\em unit-preserving $R$-matrix action} on $\Omega$ is
$$
R.\Omega = RT\Omega\ .
$$
\end{definition}

\begin{proposition}
The unit-preserving $R$-matrix action is a left group action. 
\end{proposition}

\paragraph{Proof.}
By Definition \ref{Def:RactionWithUnit},
we have
\begin{equation}\label{vwwx}
R_a. (R_b. \Omega) = R_a \bigl(z [\b1 - R_a^{-1}]\bigr) 
R_b \bigl(z [\b1 - R_b^{-1}]\bigr) \Omega.
\end{equation}
The action of $\bigl(z [\b1 - R_a^{-1}]\bigr) 
R_b$ equals the action of 
$R_b \Bigl(R_b^{-1}\bigl(z [\b1 - R_a^{-1}]\bigr) \Bigr)$ 
by Proposition~\ref{Thm:TRequalsRT}.
Thus, \eqref{vwwx} equals
\begin{multline*}
R_a R_b \Bigl(R_b^{-1}\bigl(z [\b1 - R_a^{-1}]\bigr) \Bigr)\bigl(z [\b1 - R_b^{-1}]\bigr) \Omega \\
= R_aR_b \Bigl(z \bigl[ R_b^{-1}(\b1) - R_b^{-1} R_a^{-1}(\b1) + \b1 - R_b^{-1}(\b1) \bigr] \Bigr) \Omega \\
=(R_a R_b) \bigl( z[ \b1 - (R_aR_b)^{-1}(\b1)] \bigr) \Omega\ .
\end{multline*}
Proposition \ref{bnnb} has been used in the last equality.
The result 
is precisely the definition of the unit-preserving action of~$R_aR_b$.
\qed

\section{The $R$-matrix for $A_2$} \label{Sec:A2}

We compute the $R$-matrix for the Frobenius manifold of the $A_2$ singularity and deduce an expression for the shifted $3$-spin
Witten class in terms of stable graphs. 
The outcome is a proof of Theorems 
\ref{Thm:relations} and~\ref{Thm:3Witten}.

\subsection{The Frobenius manifold $A_2$} \label{Ssec:FrobA2}

We compute all the differential geometric data associated with the Frobenius manifold $A_2$ for use in the following calculations.

The Frobenius manifold $A_2$ is based on
the 2-dimensional vector space{\footnote{In the notation
of Section \ref{wsc}, $V$ is $V_3$.}} $V$ with coordinates $x=t^0$ and $y=t^1$ corresponding to the remainders $0$ and $1$ modulo~3 respectively. The  {\em unit vector field} is $\d_x = \frac{\d}{\d x}$.
The {\em metric} is 
$$
\eta = dx \otimes dy + dy \otimes dx \ \ \ \ \ \text{or} \ \ \ \ \ 
\eta = \left( \begin{array}{cc} 0 & 1\\ 1 & 0 \end{array} \right)\ .
$$

Since the only nonzero values of  Witten's $3$-spin class in genus~$0$ are
$$
W_{0,3}(0,0,1) = 1, \qquad W_{0,4}(1,1,1,1) = \frac{1}{3}\ ,
$$
the primary genus~0 {Gromov-Witten potential} is
$$
\F(x,y) = \frac12 x^2 y + \frac1{72} y^4\ .
$$
The {\em  Euler field} is
$$
E = x \frac{\d}{\d x} + \frac23 y \frac{\d}{\d y}\ .
$$
The Lie derivatives of $E$ on the basis vectors fields are easily
calculated:
\begin{align*}
L_E(\d_x) & = [E, \d_x] = -\d_x\ ,\\
L_E(\d_y) & = [E, \d_y] = -\frac23 \d_y \ .
\end{align*}
By Proposition \ref{fvc}, 
the {\em conformal dimension} equals 
$$
\delta = \frac{r-2}{r} = \frac{1}{3}\ .
$$

Let $v$ be a tangent vector at a point of the Frobenius manifold.
We define 
the {\em shifted degree operator} $\mu(v)$, also called the {\em Hodge grading operator}, by
$$
\mu(v) = [E,v] + (1-\delta/2)v\ .
$$
Here, the vector~$v$ is extended to a flat tangent vector field in 
order to compute the commutator.
We  have
\begin{align*}
\mu(\d_x) &= -\frac16 \d_x\ ,\\
\mu(\d_y) & = \frac16 \d_y\ .
\end{align*}

\begin{definition}
To simplify the formulas, we will use the following notation:
$$
\phi = \frac{y}3\ , 
\quad \hd_x = \phi^{1/4} \d_x\ , 
\quad \hd_y = \phi^{-1/4} \d_y\ .
$$
The frame $(\hd_x, \hd_y)$ in the tangent space of~$V$ at $(x,y)$ 
is the most practical for the computations. The dual frame of the cotangent space is denoted by 
$$\hdx = \phi^{-1/4} dx\ , \ \  \hdy = \phi^{1/4} dy\ .$$
\end{definition}

The {\em quantum multiplication} of vector fields on the Frobenius manifold is given  by
\begin{align*}
\hd_x \bullet \hd_x &= \phi^{1/4} \hd_x \ ,
\\
\hd_x \bullet \hd_y &= \phi^{1/4} \hd_y\ ,
\\
\hd_y \bullet \hd_y &= \phi^{1/4} \hd_x\ .
\end{align*}

Whether in basis $(\d_x, \d_y)$ or in frame $(\hd_x, \hd_y)$, the shifted degree operator is expressed by the matrix
$$
\frac16
\begin{pmatrix}
-1 & 0\\
0 & 1
\end{pmatrix}.
$$
Unlike $\d_x$, the vector field $\hd_x$ is not flat.
However, in the definition of $\mu$, we use the flat extension of $\hd_x$ at a given point, 
which only differs from $\d_x$ by a multiplicative constant.

We will also need the operator $\xi$ of {quantum multiplication by $E$}. In the
frame $(\hd_x, \hd_y)$, $\xi$ is  given by 
$$
\xi = 
\begin{pmatrix}
\nice x & \nice 2 \phi^{3/2} \\
\\
\nice 2 \phi^{3/2} & \nice x
\end{pmatrix}.
$$

\begin{remark}
The computations not involving the Euler vector field apply more generally to 2-dimensional Frobenius manifolds whose Gromov-Witten potential has the form
$$
\F(x,y) = \frac12 x^2 y + \Phi(y)
$$
with the convention $\phi = \phi(y) = \Phi'''(y)$. For instance, the Gromov-Witten potential of ${\mathbb C} {\rm P}^1$ has the above form with 
$$\Phi = \phi = Q e^y.\ $$
\end{remark}

\subsection{The topological field theory} \label{Ssec:TopFT}

A topological field theory $\omega_{g,n}$ is a CohFT of degree~$0$, as discussed in Section~\ref{cft}. 
Teleman's reconstruction, used to prove Theorem \ref{tele}, expresses every semisimple CohFT $\Omega$ as 
a unit-preserving $R$-matrix action (see Definition~\ref{Def:RactionWithUnit}) on the 
topological field theory $\omega_{g,n}$ with unit 
where
$$\omega_{0,3} = \Omega_{0,3}.$$
Let us start by determining the topological field theory  $\omega_{g,n}$ for Witten's
$3$-spin class.

\begin{lemma}\label{nmmn}
For the topological (degree~0) part of Witten's $3$-spin theory, we have 
$$
\omega_{g,n}(\hd_x^{\otimes n_0} \otimes \hd_y^{\otimes n_1}) = 
2^g \phi^{\frac{2g-2+n}{4}} \cdot \delta^{\rm odd}_{g+n_1},
$$
where $n = n_0+n_1$.
Here, 
\begin{equation*}
\delta^{\rm odd}_{g+n_1}=
\left|
\begin{array}{cl}
1 & \mbox{ if } g+n_1 {\ is\ odd},\\
0 & \mbox{ if } g+n_1 {\ is\  even.}
\end{array}
\right.
\end{equation*}
\end{lemma}

\paragraph{Proof.} The values of $\omega_{0,3}$ are prescribed by the quantum product:
\begin{eqnarray*}
\omega_0(\hd_x \otimes \hd_x \otimes \hd_x) &=& 
\omega_0(\hd_x \otimes \hd_y \otimes \hd_y) \ \ = \ \ 0 \\
\omega_0(\hd_x \otimes \hd_x \otimes \hd_y) &=& 
\omega_0(\hd_y \otimes \hd_y \otimes \hd_y) \ \  =\ \ 
\phi^{1/4}\ .
\end{eqnarray*}
For  other $g$ and $n$, we consider a stable curve with a maximal possible number of nodes
(each component is rational with 3 special points). 
The vectors $\hd_x$ and $\hd_y$ are placed in some way on the marked points,
 and we must place either $\hd_x \otimes \hd_y$ or $\hd_y \otimes \hd_x$ at each node 
in such a way that the number of $\hd_y$'s is odd on each component of the curve. 
If $g+n_1$ is even, such a placement is impossible.
If $g+n_1$ is odd,
the placement can be done in $2^g$ ways, since the dual graph of the curve has $g$ independent cycles. 

By the factorization rules for CohFTs, the contribution of each 
successful placement of the $\hd_x$'s and $\hd_y$'s equals $\phi^{\frac{2g-2+n}{4}}$, where $2g-2+n$ is the number of rational components of the curve.
\qed

\subsection{The $R$-matrix} \label{Ssec:Rmatrix}

Givental  \cite[pages 4-5]{Givental2} gives a general
method for computing the $R$-matrix of a Frobenius manifold without using an Euler field. The method is ambiguous: 
the $R$-matrix depends on the choice of certain integration constants. 
In the presence of an Euler field~$E$, there is a unique choice of 
constants such that 
$$L_ER_m = -m R_m$$ 
for every~$m$. In the conformal case,
 Givental's method can be simplified by substituting
 $i_E$ into his recursive equation. 
The simplified method for computing the $R$-matrix of a conformal Frobenius manifold is given, for instance, by Teleman~\cite{Teleman} in the proof of the theorem of Section~8.15. Since the $3$-spin theory yields a conformal Frobenius manifold, the simplified method is suitable for us.

Let $\xi$ be the operator of quantum multiplication by the tangent vector~$E$. 
The matrices $R_m$ then 
satisfy the following recursive equation\footnote{In Teleman's paper, the commutator has the opposite sign, since his Euler field is the opposite of ours.}:
\begin{equation} \label{Eq:TelemansRecursion}
[R_{m+1},\xi] = (m+\mu) R_m.
\end{equation}
At a semisimple point of a conformal Frobenius manifold, the above
 equation determines the matrices $R_m$ uniquely starting from $R_0=1$. 
Let
$$
R_m = 
\begin{pmatrix}
a_m & b_m \\
c_m & d_m
\end{pmatrix}.
$$
Using the formulas of Section~\ref{Ssec:FrobA2} for $\xi$ and $\mu$,  we rewrite~\eqref{Eq:TelemansRecursion} as
$$
\left[
\begin{pmatrix}
a_{m+1} & b_{m+1} \\
c_{m+1} & d_{m+1}
\end{pmatrix},
\begin{pmatrix}
x & 2 \phi^{3/2}\\
2 \phi^{3/2} & x
\end{pmatrix}
\right]
=
\frac16 \begin{pmatrix}
6m-1 & 0\\
0 & 6m+1
\end{pmatrix}
\begin{pmatrix}
a_m & b_m \\
c_m & d_m
\end{pmatrix},
$$
or in other words
$$
2 \phi^{3/2}
\begin{pmatrix}
b_{m+1}-c_{m+1} & a_{m+1}-d_{m+1} \\
d_{m+1}-a_{m+1} & c_{m+1}-b_{m+1}
\end{pmatrix}
=
\frac16
\begin{pmatrix}
(6m-1) a_m & (6m-1) b_m \\
(6m+1) c_m & (6m+1) d_m
\end{pmatrix}.
$$
The following formulas are easily checked to
be the unique solutions:
\begin{align*}
a_m &= \frac1{1728^m\, \phi^{3m/2}} \; 
\frac{1+6m}{1-6m} \; \frac{(6m)!}{(3m)!\,(2m)!} \; \de,\\
b_m &= \frac1{1728^m\, \phi^{3m/2}} \;  
\frac{1+6m}{1-6m} \; \frac{(6m)!}{(3m)!\,(2m)!} \; \dod,\\
c_m &= \frac1{1728^m\, \phi^{3m/2}} \;  
\frac{(6m)!}{(3m)!\,(2m)!} \; \dod,\\
d_m &= \frac1{1728^m\, \phi^{3m/2}} \;  
\frac{(6m)!}{(3m)!\,(2m)!} \; \de.
\end{align*}

We now make explicit the connection with the central
power series discovered by Faber and Zagier,
$$
\B_0(T) = \sum_{m \geq 0} \frac{(6m)!}{(2m)!(3m)!}(-T)^m,
\ \
\B_1(T) = \sum_{m \geq 0} \frac{1+6m}{1-6m} \frac{(6m)!}{(2m)!(3m)!}(-T)^m.
$$
Denote by $\Be_0$, $\Bo_0$, $\Be_1$, and $\Bo_1$ the respective
 even and odd degree parts.
The final expression for the $R$-matrix is:
\begin{equation}\label{gvcc}
R(z) = 
\begin{pmatrix}
\nice \Be_1 \left( \frac{z}{1728 \, \phi^{3/2}} \right) &
\nice -\Bo_1 \left( \frac{z}{1728 \, \phi^{3/2}} \right) \\
\\
\nice -\Bo_0 \left( \frac{z}{1728 \, \phi^{3/2}} \right) &
\nice \Be_0 \left( \frac{z}{1728 \, \phi^{3/2}} \right)
\end{pmatrix}.
\end{equation}
The symplectic condition for the $R$-matrix follows from the identity
$$
\B_0(T) \B_1(-T) + \B_0(-T) \B_1(T) =2,
$$
or, equivalently,
$$
\Be_0(T) \Be_1(T) - \Bo_0(T) \Bo_1(T) = 1
$$
discovered previously in~\cite{Pixton}.
Using the identity, we find
\begin{equation}\label{vssa}
R^{-1}(z) = 
\begin{pmatrix}
\nice \Be_0 \left( \frac{z}{1728 \, \phi^{3/2}} \right) &
\nice \Bo_1 \left( \frac{z}{1728 \, \phi^{3/2}} \right) \\
\\
\nice \Bo_0 \left( \frac{z}{1728 \, \phi^{3/2}} \right) &
\nice \Be_1 \left( \frac{z}{1728 \, \phi^{3/2}} \right)
\end{pmatrix}.
\end{equation}

\subsection{An expression for the shifted $3$-spin Witten class}
\label{Ssec:expression}

We combine here the expression for the topological field theory from Section~\ref{Ssec:TopFT} with the $R$-matrix action from Definition~\ref{Def:RactionWithUnit} using the explicit formulas for the $R$-matrix of Section 
\ref{Ssec:Rmatrix}.

Let $\tau = (x,y)$, $y \ne 0$, be a point of the Frobenius manifold~$A_2$. 
Let $a_1, \dots, a_n \in \{0,1 \}$ and let 
$$D = \frac{g-1 + \sum_{i=1}^n a_i}{3}$$
 be the degree of Witten's $3$-spin class. 
By convention, $\phi = y/3$. 
Recall the expressions $\cR^d_{g, (a_1, \dots, a_n)}$ of
Definition~\ref{Not:relations}.

\begin{theorem} \label{Thm:expression}
Witten's class for the shifted $3$-spin theory equals
$$
\W_{g,n}^\tau(\d_{a_1} \otimes \cdots \otimes \d_{a_n})= 2^g
\sum_{d \geq 0} \frac{\phi^{\frac32(D-d)}}{1728^d} \;
q\left(\cR^d_{g,(a_1, \dots, a_n)}\right),
$$
where $\d_0 = \d_x$, $\d_1 = \d_y$.
\end{theorem}

The following Corollary is an immediate consequence of Theorem 
\ref{Thm:expression} and the equation
$$
\W_{g,n}^\tau(\d_{a_1} \otimes \cdots \otimes \d_{a_n}) = W_{g,n}(a_1, \dots, a_n)
+ \mbox{ lower degree terms.}
$$
explained in Section \ref{Sec:FrobA}.
Theorems 
\ref{Thm:relations} and~\ref{Thm:3Witten}
are implied by the Corollary.

\begin{corollary}
We have the evaluations:
\begin{align*}
q\left(\cR^d_{g,(a_1, \dots, a_n)}\right) 
& = 2^g \, 1728^D \, W_{g,n}(a_1, \dots, a_n) 
&  \mbox{\rm for} \quad d=D,\\
q\left(\cR^d_{g,(a_1, \dots, a_n)}\right) & = 0 & 
 \mbox{\rm for} \quad d > D.
\end{align*}
\end{corollary}

\paragraph{Proof of Theorem~\ref{Thm:expression}.}
By Teleman's reconstruction result
in the conformal semisimple case,
Witten's shifted $3$-spin class is given by $R.\omega$
where
\begin{enumerate}
\item[$\bullet$] $R$ is given by \eqref{gvcc},
\item[$\bullet$] $\omega$ is the topological part of the
shifted $3$-spin theory.
\end{enumerate} 
The proof now just amounts to a systematic matching of all factors 
in the sums over stable graphs which occur in Definition
\ref{Def:RactionWithUnit} for the $R$-matrix action and Definition 
\ref{Not:relations} 
for $\cR^d_{g,(a_1, \dots, a_n)}$. 

Consider first the expression for 
the CohFT $R.\omega$
applied to a tensor product of $n$ vectors $\d_x$ and $\d_y$. 
As before, we denote by $n_0$ and $n_1$ the number of $0$s and $1$s among $a_1, \dots, a_n$ so that $n_0+n_1=n$. 

\vspace{8pt}
\noindent{\bf Powers of $\phi$.} 
\vspace{8pt}

Since we wrote the $R$-matrix in frame $(\hd_x, \hd_y)$, 
we must substitute
 $$\d_x \mapsto \phi^{-1/4} \hd_x\ , \ \ \ 
\d_y \mapsto \phi^{1/4} \hd_y\ $$
in the tensor product argument for $R.\omega$.
 The result of the substitution is a factor of $\phi^{\frac{n_1-n_0}{4}}$. 

By formula~\eqref{vssa},
all coefficients of $R^{-1}_m$ contain a factor of $\phi^{-3m/2}$.
Tracing through the definitions of the all the actions 
\begin{equation}\label{gzzq}
R.\omega = R T \omega\ , \ \ \ T(z) = z\cdot[\partial_x - 
R^{-1}(\partial_x)](z),
\end{equation}
the $R$-matrix contributes a factor of~$\phi^{-3d/2}$, where $d$ is the degree 
of the class.

By Lemma \ref{nmmn},
the topological field theory $\omega$ contributes (subject to parity condition
accounted for later) a factor of 
$\phi^{\frac{2\mathsf{g}_v-2+\mathsf{n}_v}{4}}$ for every vertex~$v$. 
These factors combine to yield $\phi^{\frac{2g-2+n}{4}}$.

Finally, each $\kappa$-leg contributes in two way.
First, since we must substitute
$$\d_x \mapsto \phi^{-1/4} \hd_x$$
in formula  \eqref{gzzq} for $T(z)$, each $\kappa$-leg contributes
 $\phi^{-1/4}$.
Second, because the $\kappa$-leg increases the valence of the vertex by 1,
a factor of  
$\phi^{1/4}$ is contributed via the topological field theory. 
Thus, the contributions of each $\kappa$-leg to the power of $\phi$
cancel. 

Collecting all of the above factors, we obtain
a final calculation of the exponent of $\phi$:
$$
\frac{n_1-n_0}4 -\frac{3d}2 + \frac{2g-2+n}4
=
\frac{g-1+n_1-3d}2 = \frac{3D-3d}2 = \frac32(D-d).
$$

\noindent {\bf Powers of 1728.}
\vspace{8pt} 

All coefficients of $R^{-1}_m$ contain a factor of $1/1728$. 
Hence, as above, we obtain a factor of~$1728^{-d}$ from the
$R$-matrix action.

\vspace{8pt}
\noindent {\bf Powers of~2.} 
\vspace{8pt}

At each vertex the topological field theory contributes a factor of $2^{g_v}$. These combine into
$$
\prod_{v\in V(\Gamma)} 2^{g_v} = \frac{2^g}{2^{h^1(\Gamma)}}.
$$
The factor $2^{-h^1(\Gamma)}$ is present in the definition of $\cR^d_{g,(a_1, \dots, a_n)}$, and the remaining $2^g$ is included in the
statement of Theorem~\ref{Thm:expression}.

\vspace{8pt}
\noindent{\bf Parity conditions at the vertices.} 
\vspace{8pt}

The topological field theory $\omega$
provides a nonzero contribution at a vertex if and only if $g_v+n_1(v)$ is odd. We must  prove
the parity condition which occurs in the definition of $\cR^d_{g, (a_1, \dots, a_n)}$ exactly
matches.

The parity condition is imposed on $\cR^d_{g,(a_1, \dots, a_n)}$ by extracting the coefficient of $\zeta_v^{g_v-1}$,
at each vertex $v$: see Definition \ref{Not:relations}.
We may view the factors of 
$\zeta_v$ as having the following sources. 
A leg carrying the assignment $a_l=1$ (corresponding to $\d_y$)
 contributes a $\zeta_v$, while a leg carrying the assignment $a_l=0$ (corresponding 
to $\d_x$) does not. 
The terms of $\Bo_0$ (including the effect of the $\kappa$-legs) and the terms of $\Bo_1$ contribute
 a $\zeta_v$. The terms of $\Be_0$ and $\Be_1$ do not contribute anything (because they leave the parity invariant). 
Finally, every edge insertion $\Delta_e$ contributes a factor if $e$ is adjacent to~$v$. 
The edge term of Definition \ref{Not:relations}  can be expanded 
via  
$$\B_0=\Be_0+\Bo_0 \ , \ \ \ \B_1=\Be_1+\Bo_1\ $$
and matched with the edge term of the CohFT $R.\omega$ using \eqref{Eq:zw} and \eqref{vssa}.
Then the contributing factor is
 $\zeta_v$ if the bi-vector includes a factor $\hd_y$ on the side of the vertex~$v$ and 1 otherwise. 
Hence, the power of the variable $\zeta_v$ correctly counts the parity of entries $\hd_y$ submitted
to the topological field theory $\omega$ at the vertex~$v$.

\vspace{8pt}
\noindent{\bf Coefficients of the series $\B$.} 
\vspace{8pt}

These coefficients simply coincide in the expression for $\cR^d_{g,(a_1, \dots, a_n)}$ and the formulas of the unit-preserving $R$-matrix action in all instances (legs, $\kappa$-legs, and edges).
 \qed

\subsection{$\tP$ implies $\P$}\label{section:cor2}

We present here the proof of Corollary \ref{Pix}: the derivation of the more complete
set of relations $\P$ conjectured in \cite{Pixton} from the 
set 
$\tP$ proven in Theorem \ref{Thm:relations}.

Our relations $\cR^d_{g,A}$ differ from the relations $\cR^d_{g,A,\sigma}$ of~\cite{Pixton} in three ways. 
First, the signs of the coefficients in the series $\B_0$ and $\B_1$ are modified. 
The outcome is a global change of sign in some of the relations. 
Second, the range of the $a_i$'s is different. In our relations, the $a_i$'s are equal to 0 or 1, while in~\cite{Pixton},  the $a_i$'s can be any integers equal to 0 or 1 modulo 3. In fact, replacing an $a_i$ by $a_i+3$ in the relations of $\P$ amounts to multiplying the relation by~$\psi_i$. 
Therefore taking $a_i < 3$ is sufficient. 
Finally, the relations of~\cite{Pixton} also depend on a partition $\sigma$. In our relations, we are 
implicitly considering only the empty partition case. 
A relation with a nonempty partition $\sigma$ is easily
 obtained from a relation with an empty $\sigma$ by push-forward:
$$
\cR^d_{g,(a_1, \dots, a_n), (\sigma_1, \dots, \sigma_m)} = 
p_* \cR^d_{g,(a_1, \dots, a_n, \sigma_1+3, \dots, \sigma_m+3)},
$$
where $p:\oM_{g,n+m} \to \oM_{g,n}$ is the forgetful map.
Thus the relations $\cR^d_{g,A}$ imply all the relations $\cR^d_{g,A,\sigma}$. 

The span of $\tP$ is not an ideal, but generates an ideal in each~$\cS_{g,n}$. 
The associated family of ideals is closed under pull-backs by forgetful maps 
(because of axiom~(iii) of a CohFT with unit) and gluing maps (because of axiom~(ii) of a CohFT). The family of ideals
 is not closed under push-forwards by forgetful maps and gluing maps. 
After taking the closure under push-forwards by forgetful maps,
 we obtain the span of $\cR^d_{g,A,\sigma}$, as we have just proved. Taking the closure under push-forwards by gluing maps we get the full set of relations~$\P$ from~\cite{Pixton}. \qed

\subsection{Examples} \label{Ssec:examples}
\begin{example}
Let $g=0$, $n=3$. Here, we have
$$
\W_{0,3}^\tau(\d_x \otimes \d_x \otimes \d_y) = 1, \qquad
\W_{0,3}^\tau(\d_y \otimes \d_y \otimes \d_y) = \phi = \frac{y}3.
$$
These values come directly from the topological field theory -- the $R$-matrix is not needed.
The first expression equals the Witten class $W_{0,3}(0,0,1)$, and the second expression is the push-forward of 
$$
\W_{0,4}(\d_y \otimes \d_y \otimes \d_y \otimes y \d_y) = yW_{0,4}(1,1,1,1).
$$
In both cases, no further $y\d_y$ insertions are possible for dimension reasons.
\end{example}

\begin{example}
Let $g=0$, $n=4$. We will study all 5 cases.

\medskip

{\bf First case:} $\W_{0,4}^\tau(\d_x \otimes \d_x \otimes \d_x \otimes \d_x)$. We have $D = -1/3$. The parity condition imposes $d=1$. Thus $\frac32(D-d) = -2$. By Theorem \ref{Thm:expression}, we find
\begin{equation*} 
\W_{0,4}^\tau(\d_x \otimes \d_x \otimes \d_x \otimes \d_x)
= \frac{60 \kappa_1 - 60 \sum_{i=1}^4 \psi_i +60 \delta}
{1728 \phi^2}.
\end{equation*}
Since $d>D$, the above expression  must be~0.
We obtain the first nontrivial relation:
$$
\kappa_1 -\sum_{i=1}^4\psi_i + \delta = 0 \ \in H^2(\oM_{0,4}, \mathbb{Q})\ .$$
The relation is true by the following basic evaluation in $H^2(\oM_{0,4}, \mathbb{Q})$:
$$\kappa_1 = [{\rm point}],\ \ \psi_i = [{\rm point}], \ \ \delta = 3 [{\rm point}] \ .$$
Alternatively, we see that this expression coincides up to a factor with Mumford's formula for $\lambda_1$,
 and  $\lambda_1=0$ in genus~0.

\medskip

{\bf Second case:} $\W_{0,4}^\tau(\d_x \otimes \d_x \otimes \d_x \otimes \d_y)$. We have $D = 0$. The parity condition imposes $d=0$. Thus $\frac32(D-d) = 0$. 
From Theorem \ref{Thm:expression}, we obtain
$$
\W_{0,4}^\tau(\d_x \otimes \d_x \otimes \d_x \otimes \d_y)
= 1.
$$
Since $d=D$ we know that this expression should be equal to Witten's class, which is indeed the case: $W_{0,4}(0,0,0,1) = 1$.

\medskip

{\bf Third case:} $\W_{0,4}^\tau(\d_x \otimes \d_x \otimes \d_y \otimes \d_y)$. We have $D = 1/3$. The parity condition imposes $d=1$. Thus $\frac32(D-d) = -1$. We obtain
$$
\W_{0,4}^\tau(\d_x \otimes \d_x \otimes \d_y \otimes \d_y)=
$$
$$
\frac{60 \kappa_1 - 60 (\psi_1+\psi_2) + 84(\psi_3+\psi_4) +60 \delta_{[1,2|3,4]} - 84 (\delta_{[1,3|2,4]} + \delta_{[1,4|2,3]})}
{1728 \phi}.
$$
Since $d>D$, the expression must vanish (as is easly checked).

\medskip

{\bf Fourth case:} $\W_{0,4}^\tau(\d_x \otimes \d_y \otimes \d_y \otimes \d_y)$. We have $D = 2/3$. The parity condition imposes $d=0$. Thus $\frac32(D-d) = 1$. We obtain
$$
\W_{0,4}^\tau(\d_x \otimes \d_y \otimes \d_y \otimes \d_y)= \phi
$$
which is the push-forward of 
$$
\W_{0,5}(\d_x \otimes \d_y \otimes \d_y \otimes \d_y \otimes y \d_y) 
=y W_{0,5}(0,1,1,1,1)
$$ 
under the forgetful map forgetting the last marked point.

\medskip

{\bf Fifth case:} $\W_{0,4}^\tau(\d_y \otimes \d_y \otimes \d_y \otimes \d_y)$. We have $D = 1$. The parity condition imposes $d=1$. Thus $\frac32(D-d) = 0$. We obtain
\begin{eqnarray*}
\W_{0,4}^\tau(\d_y \otimes \d_y \otimes \d_y \otimes \d_y) & = & 
\frac{60 \kappa_1 + 84 \sum_{i=1}^4 \psi_i + 60 \delta}{1728} \\
& =& \frac{(60 + 84 \cdot 4 + 60 \cdot 3) [\mbox{pt}]}{1728} \\
& = & \frac13 [\mbox{point}]
\end{eqnarray*}
which is the correct value of Witten's class 
$$
W_{0,4}(1,1,1,1)=\frac13  [\mbox{point}].
$$

The relations obtained through these computations in $d=1$ (after dividing by 12 or by 60)
are listed below:
\begin{align*}
\kappa_1 \; -\; \psi_1 \;-\; \psi_2 \;-\; \psi_3 \;-\; \psi_4 \;+\; \delta_{[1,2|3,4]} \;+\; \delta_{[1,3|2,4]} \;+\; \delta_{[1,4|2,3]} & = 0, \\
5\kappa_1 - 5\psi_1 - 5\psi_2 + 7\psi_3 + 7\psi_4 + 5\delta_{[1,2|3,4]} -7\delta_{[1,3|2,4]} -7 \delta_{[1,4|2,3]} & = 0, \\
5\kappa_1 +7 \psi_1 + 7 \psi_2 - 5\psi_3 - 5 \psi_4 + 5\delta_{[1,2|3,4]} -7\delta_{[1,3|2,4]} -7 \delta_{[1,4|2,3]} & = 0, \\
5\kappa_1 - 5\psi_1 +7 \psi_2 -5 \psi_3 + 7\psi_4 -7 \delta_{[1,2|3,4]} +5 \delta_{[1,3|2,4]} -7 \delta_{[1,4|2,3]} & = 0, \\
5\kappa_1 +7\psi_1 -5 \psi_2 +7 \psi_3 -5 \psi_4 -7 \delta_{[1,2|3,4]} +5 \delta_{[1,3|2,4]} -7 \delta_{[1,4|2,3]} & = 0, \\
5\kappa_1 - 5\psi_1 +7 \psi_2 +7 \psi_3 -5 \psi_4 -7 \delta_{[1,2|3,4]} -7 \delta_{[1,3|2,4]} +5 \delta_{[1,4|2,3]} & = 0, \\
5\kappa_1 +7\psi_1 -5 \psi_2 -5 \psi_3 +7 \psi_4 -7 \delta_{[1,2|3,4]} -7 \delta_{[1,3|2,4]} +5 \delta_{[1,4|2,3]} & = 0.\\
\end{align*}
After some linear algebra, the system is equivalent to:
$$
\kappa_1 = \psi_1 = \psi_2 = \psi_3 = \psi_4 = \delta_{[1,2|3,4]} = \delta_{[1,3|2,4]} = \delta_{[1,4|2,3]}.
$$
We have obtained a complete set of relations in $RH^2(\oM_{0,4})$.
\end{example}

\begin{example}
The Getzler relation \cite{Get} is a degree~2 relation in $\cS_{1,4}$ 
which can not be obtained by the pull-back of any simpler relations. 
Since
$$ 2 > \frac{1-1+ 1+1+1+1}{3} = \frac{4}{3}\ ,$$
the relation
$\cR^2_{1,(1,1,1,1)}$ lies in the set $\tP$. In fact, $\cR^2_{1,(1,1,1,1)}$
is the Getzler relation (modulo more elementary genus 0 and 1 relations).
 
The Belorousski-Pandharipande relation \cite{bp}
is a degree~2 relation in $\cS_{2,3}$ which can not be obtained by the
pull-back of any simpler relations. 
The relation $\cR^2_{2,(1,1,1)}$ lies in $\tP$ since
$$ 2 > \frac{2-1+ 1+1+1}{3} = \frac{4}{3}\ $$
and 
is an equivalent form of the BP equation.
\end{example}

The outcome of several such investigations is reported in \cite{Pixton}. All known relations have been
explained by Theorem 
\ref{Thm:expression}.

\subsection{Some concluding remarks} \label{Ssec:remarks}

Our computations provide an instructive example of what happens to a Coh\-FT as we move towards a non-semisimple point of a Frobenius manifold. Let us examine more closely the limit of our expressions for the 
shifted Witten $3$-spin class as $y \to 0$ or, in other words, $\phi \to 0$. 
The coefficients of the $R$-matrix involve negative powers of $\phi$, therefore the
$R$-matrix diverges. 
The topological field theory $\omega$ to which we apply the $R$-matrix involves positive powers of $\phi$. 
As a result, each term of our expression for the shifted Witten class comes with a factor
$$
\phi^{\frac23 (D - d)},
$$
where $D = \D_{g,n}(a_1, \dots, a_n)$ is the degree of Witten's class and $d$ is the degree of the term in question. As $\phi \to 0$, the terms of degree less then $D$ tend to~0,
 the terms of degree equal to~$D$ are invariant, and  
the terms of degree greater than~$D$ diverge. At first sight the expression appears
 to diverge, but {\em because the terms of degree greater than~$D$ combine into tautological relations} the expression actually has a finite limit equal to Witten's $3$-spin class.

A natural question is 
 whether our formulas for Witten's class lift from $\oM_{g,n}$ to formulas on 
the space $\oM^{1/r}_{g; a_1, \dots, a_n}$ of $r$-spin structures. The answer is {\em no}: 
the divisibility condition 
$$\frac{(g-1)(r-2)+ \sum_i a_i}{r} \ \in {\mathbb Z}$$ 
does not necessarily hold for each vertex of the dual graph.
Hence, there is no natural boundary stratum in $\oM^{1/r}_{g; a_1, \dots, a_n}$ where the terms of our formula can be lifted. 
Moreover, in the simplest case $r=2$, we have
$$
\W_{g,n}(0, \dots, 0) = 
\left|
\begin{array}{rl}
1 & \mbox{ if the spin structure is even,} \\
-1 &  \mbox{ if the spin structure is odd.}
\end{array}
\right.
$$
Such an  answer cannot be expressed in terms of dual graphs at all. Some more structure is
required.

\vspace{+16 pt}
\noindent Departement Mathematik \\
\noindent ETH Z\"urich \\
\noindent rahul@math.ethz.ch

\vspace{+8 pt}
\noindent
Department of Mathematics\\
Princeton University\\
apixton@math.princeton.edu

\vspace{+8 pt}
\noindent
CNRS, Institut Math\'ematique de Jussieu\\
zvonkine@math.jussieu.fr

\end{document}